\documentclass[11pt,a4paper]{article}
\usepackage{amsmath}
\usepackage{amsthm}
\usepackage{amssymb}
\usepackage{dsfont}
\usepackage[T1]{fontenc}
\usepackage[french,english,russian]{babel}
\usepackage[utf8]{inputenc}
\usepackage{mathrsfs}
\usepackage{lmodern}
\usepackage{graphicx}
\newtheorem{lemma}{Lemma}

\newtheorem{prop}{Proposition}
\newtheorem{theorem}{Theorem}
\newtheorem*{theorems}{Theorem}

\newtheorem{defi}{Definition}
\makeatletter
\def\moverlay{\mathpalette\mov@rlay}
\def\mov@rlay#1#2{\leavevmode\vtop{%
   \baselineskip\z@skip \lineskiplimit-\maxdimen
   \ialign{\hfil$\m@th#1##$\hfil\cr#2\crcr}}}
\newcommand{\charfusion}[3][\mathord]{
    #1{\ifx#1\mathop\vphantom{#2}\fi
        \mathpalette\mov@rlay{#2\cr#3}
      }
    \ifx#1\mathop\expandafter\displaylimits\fi}
\makeatother

\title{Limit laws in the lattice problem. \\
III. Return to the case of boxes}
\author{Julien Trevisan}
\begin{document}
\selectlanguage{english}
\maketitle
\bigskip
\section*{}
\selectlanguage{french}
(Résumé en français)Nous étudions l'erreur du nombre de points d'un réseau $L$ qui appartiennent à un rectangle, centré en $0$, dont les axes sont parallèles aux axes de coordonnées, dilaté d'un facteur $t$ puis translaté d'un vecteur $X \in \mathbb{R}^{2}$. Quand nous considérons le moment d'ordre 2 de l'erreur relativement à $X \in \mathbb{R}^{2}/L$, on montre que, quand $t$ est aléatoire et devient grand et quand l'erreur est normalisée par une quantité qui se comporte, dans le cas admissible, comme $\sqrt{\log(t)}$,  elle converge en loi vers une constante positive explicite. Dans le cas d'un réseau $L$ typique, on montre que ce résultat tient toujours mais la normalisation est toutefois plus importante, autour de $\log(t)$. On montre aussi que quand $L = \mathbb{Z}^{2}$, l'erreur, quand elle est normalisée par $t$, converge en loi quand $t$ est aléatoire et devient grand et on calcule les moments de la loi limite. \\
\\
\selectlanguage{english}
(Engish abstract)We study the error of the number of points of a lattice $L$ that belong to a rectangle, centred at $0$, whose axes are parallel to the coordinate axes, dilated by a factor $t$ and then translated by a vector $X \in \mathbb{R}^{2}$. When we consider the second order moment of the error relatively to $X \in \mathbb{R}^{2}/L$, one shows that, when $t$ is random and becomes large and when the error is normalized by a quantity which behaves, in the admissible case, as $\sqrt{\log(t)}$, it converges in distribution to an explicit positive constant. In the case of a typical lattice $L$, we show that this result still holds but the normalisation is more important, around $\log(t)$. We also show that when $L=\mathbb{Z}^{2}$, the error, when normalized by $t$, converges in distribution when $t$ is random and becomes large and we compute the moments of the limit distribution.
 \selectlanguage{english} 
\section{Introduction}
 Let $P$ be a measurable subset of $\mathbb{R}^{d}$ of non-zero finite Lebesgue measure. We want to evaluate the following cardinal number when $t \rightarrow \infty$: $$ N(tP + X, L) = | (tP + X) \cap L|$$ where $X \in \mathbb{R}^{d}$, $L$ is a lattice of $\mathbb{R}^{d}$ and $t P + X$ denotes the set $P$ dilated by a factor $t$ relatively to $0$ and then translated by the vector $X$.  \\
Under mild regularity conditions on the set $P$, one can show that: 
$$ N(tP + X, L) = t^{d}\frac{\text{Vol}(P)}{\text{Covol}(L)} + o(t^{d}) $$
where $o(f(t))$ denotes a quantity such that, when divided by $f(t)$, it goes to $0$ when $t \rightarrow \infty$ and where $\text{Covol}(L)$ is the volume of a fundamental set of the lattice $L$. \\
We are interested in the error term $$\mathcal{R}(tP + X,L) = N(tP + X, L) - t^{d}\frac{\text{Vol}(P)}{\text{Covol}(L)} \textit{.}$$
In the case where $d=2$ and where $P$ is the unit disk $\mathbb{D}^{2}$, Hardy's conjecture in $\cite{hardy1917average}$ stipulates that we should have for all $\epsilon > 0$, $$\mathcal{R}(t \mathbb{D}^{2}, \mathbb{Z}^{2}) = O(t^{\frac{1}{2}+\epsilon}) $$
where $Y = O(X)$ means that there exists $D > 0$ such that $ |Y| \leqslant D |X|$. \\
One of the result in this direction has been established by Iwaniec and Mozzochi in $\cite{iwaniec1988divisor}$. They have proven that for all $\epsilon > 0$, $$\mathcal{R}(t \mathbb{D}^{2}, \mathbb{Z}^{2}) = O(t^{\frac{7}{11}+\epsilon}) \textit{.}$$
This result has been recently improved by Huxley in \cite{huxley2003exponential}. Indeed, he has proven that: $$\mathcal{R}(t \mathbb{D}^{2}, \mathbb{Z}^{2}) = O(t^{K} \log(t)^{\Lambda}) $$
where $K = \frac{131}{208} $ and $\Lambda = \frac{18627}{8320}$.  \\
In dimension 3, Heath-Brown has proven in $\cite{heath2012lattice}$ that: 
$$\mathcal{R}(t \mathbb{D}^{3}, \mathbb{Z}^{3}) = O(t^{\frac{21}{16}+\epsilon}) \textit{.}$$
These last two results are all based on estimating what are called $\textit{exponential sums}$.\\ Furthermore, in both cases, the error is considered in a deterministic way. \\
Another approach was followed first by Heath-Brown in $\cite{heath1992distribution}$ and then by Bleher, Cheng, Dyson and Lebowitz in $\cite{bleher1993distribution}$. They took interest in the case where the dilatation parameter $t$ is random. More precisely, they assumed that $t$ was being distributed according to the measure $\rho(\frac{t}{T}) dt$ (that is absolutely continuous relatively to Lebesgue measure) and where $\rho$ is a probability density on $[0,1]$ and $T$ is parameter that goes to infinity. In that case, Bleher, Cheng, Dyson and Lebowitz showed the following result (which generalizes the result of Heath-Brown): 
\begin{theorems}[$\cite{bleher1993distribution}$]
Let $\alpha \in [0,1[^{2}$. There exists a probability density $p_{\alpha}$ on $\mathbb{R}$ such that for every piecewise continuous and bounded function $g: \mathbb{R} \longrightarrow \mathbb{R}$, 
$$ \lim_{T \rightarrow \infty} \frac{1}{T} \int_{0}^{T} g \big( \frac{\mathcal{R}(t \mathbb{D}^{2} + \alpha , \mathbb{Z}^{2})}{\sqrt{t}} \big) \rho(\frac{t}{T}) dt = \int_{\mathbb{R}} g(x) p_{\alpha}(x) dx \textit{.}$$
Furthermore $p_{\alpha}$ can be extended as an analytic function over $\mathbb{C}$ and satisfies that for every $\epsilon > 0$, $$p_{\alpha}(x) = O(e^{-|x|^{4- \epsilon}})$$ when $x \in \mathbb{R}$ and when $|x| \rightarrow \infty$. 
\end{theorems}
We want to follow this approach on another problem. Namely, let us take $a > 0$ and $b > 0$ and let us define $\text{Rect}(a,b)$ the rectangle centred around $(0,0)$ whose summits are $(a,b)$, $(-a,b)$, $(-a,-b)$ and $(a,-b)$. \\
Let us recall the following definitions: 
\begin{defi}
\label{chap4:def2}
For a lattice $L$ of $\mathbb{R}^{d}$, its dual lattice $L^{\perp}$ is defined by 
$$L^{\perp} = \{ x \in \mathbb{R}^{d} \text { } | \text{ } \forall l \in L \text{, } <l,x> \in \mathbb{Z} \} $$
where $< , >$ is the usual euclidean scalar product over $\mathbb{R}^{d}$.
\end{defi}

\begin{defi}
\label{chap4:def1}
A lattice $L$ of $\mathbb{R}^{d}$ is called admissible if there exists $C > 0$ such that for all $l=(l_{1},\cdots,l_{d}) \in L-\{ 0 \}$, $$| \text{Num}(l)|  \geqslant C  $$ 
where $\text{Num}(l) = l_{1} \cdots l_{d}$. \\
During the rest of this article, we are going to use the following notation: 
$$\text{Num}(L) = \inf \{ | \text{Num}(l) | \text{ } | \text{ }  l \in L-\{ 0 \} \} \textit{.} $$
With this notation, saying that $L$ is admissible is equivalent to saying that $\text{Num}(L) > 0$. 
\end{defi}
We recall that if $L$ is an admissible lattice, $L^{\perp}$ is also admissible (see, for example, $\cite{skriganov1994constructions}$). \\
Let us also define the following quantity:
\begin{defi}
\label{chap4:def10}
For every lattice $L$ of $\mathbb{R}^{2}$, for every $t > 0$, one sets: 
$$V(L,t) = \sum_{\substack{ l \in L \\ 0 < \lVert l \rVert \leqslant t }} \frac{1}{l_{1}^{2} l_{2}^{2}}  $$
when for all $l \in L$ such that $ 0 < \lVert l \rVert \leqslant t $, $\text{Num}(l) \neq 0$ and where $\lVert \cdot \rVert$ is the usual euclidean norm over $\mathbb{R}^{2}$.
\end{defi}
We are going to use the following notation: for $f(t,X)$ a function of $\mathbb{R}^{+} \times \mathbb{R}^{2}$ such that $f(t,\cdot)$ is $L$-periodic, 
$$\mathbb{E}_{X \in \mathbb{R}^{2}/L}(f)(t) = \frac{1}{\textit{covol}(L)} \int_{X \in \mathbb{R}^{2} / L } f(t,X) dX \textit{.} $$ 
Let us take a probability density $\rho$ over $[0,1]$. \\
One of the goal of this article is to prove the following theorem: 
\begin{theorem}
\label{chap4:thm1} %peut-être mettre en fonction de T en fait après
Let $L$ be an admissible lattice of $\mathbb{R}^{2}$. \\
Then, first, there exists $C > 0$ such that for all $t$ large enough: $$ \frac{1}{C} \log(t) \leqslant V(L,t) \leqslant C \log(t) $$ and we are going to write it $V(L,t) = \Theta(\log(t))$. \\
Second, if $P=\text{Rect}(a,b)$, 
$$ \mathbb{E}_{X \in \mathbb{R}^{2}/L}\left( (\frac{\mathcal{R}(t P + X ,L)}{\sqrt{V(L^{\perp},t)}})^{2} \right) \textit{ converges in distribution and in probability towards } \frac{1}{4 \pi^{4} \text{Covol}(L)^{2}} $$
when $t \in [0, T]$ is distributed according to $\frac{1}{T} \rho(\frac{t}{T}) dt$ and when $T \rightarrow \infty$.
\end{theorem}
Theorem $\ref{chap4:thm1}$ can be interpreted as followed: after averaging on $X$, when $t$ is large and random, the error committed by making the approximation $$N(t \text{Rect}(a,b) + X, \Gamma_{\alpha,\alpha'}) \approx t^{2} \frac{\text{Area}(\text{Rect}(a,b))}{\text{Covol}(\Gamma_{\alpha,\alpha'})}$$ is of order $\sqrt{\log(t)}$. This normalization is in fact suggested by $\cite{trevisan2021}$. \\
This theorem is also in fact suggested by Theorem $\ref{chap4:thm2}$ that we will prove (see the next section for some heuristic explanations). We are going to state it. \\ Let us designate the space of unimodular lattices of $\mathbb{R}^{d}$ by the notation $\mathscr{S}_{d}$. \\
We are calling $\Delta_{r}$ the following set: 
$$\Delta_{r} = \{  \text{Diag}(e^{t_{1}}, \cdots, e^{t_{d}}) \text{ } | \text{ } (t_{1},\cdots,t_{d}) \in \mathbb{Z}^{d} \text{, } t_{1}+ \cdots + t_{d} = 0 \textit{ and } \lVert (t_{1},\cdots,t_{d}) \rVert \leqslant r \} \textit{.} $$ 
We have that $\Delta_{r} \subset \mathscr{S}_{d}$ and $|\Delta_{r}| = n_{r}r^{d-1} + o(r^{d-1})$ when $r \rightarrow \infty$. \\
We also call 
\begin{equation}
\label{chap4:ajout}
\Delta =  \{  \text{Diag}(e^{t_{1}}, \cdots, e^{t_{d}}) \text{ } | \text{ } (t_{1},\cdots,t_{d}) \in \mathbb{Z}^{d} \text{ and } t_{1}+ \cdots + t_{d} = 0  \}\textit{ and}
\end{equation}
for every lattice $L$ of $\mathbb{R}^{d}$, we call 
$$\lVert L \rVert = \min \{ \lVert l \rVert \text{ } | \text{ } l \in L-\{ 0 \} \} \textit{.} $$
With these notations, the Theorem $\ref{chap4:thm2}$ states as the following: 
\begin{theorem}
\label{chap4:thm2}
Let $L \in \mathscr{S}_{d}$ be an admissible lattice. Let us take $(\Theta_{\delta})_{\delta \in \Delta}$ a sequence of independent identically distributed real random variables that are symmetrical and admit moment of order $3$ and whose $\textit{alea}$ $\omega$ belongs to a probability space $\Omega$. \\
Let us set:  $$ \tilde{V}(L,r) = \sum_{\delta \in \Delta_{r}} \frac{1}{\lVert \delta L \rVert^{2 d}} $$ 
and 
$$\tilde{S}(L,\omega,r) =  \sum_{\delta \in \Delta_{r}} \frac{\theta_{\delta}(\omega)}{\lVert \delta L \rVert^{d}} \textit{.}$$
 Then one has: $$ \tilde{V}(L,r) = \Theta(r^{d-1}) $$ 
and $ \frac{\tilde{S}(L,\omega,r)}{\sqrt{\tilde{V}(L,r)}}$ converges in distribution towards the standard normal distribution when $r \rightarrow \infty$. 
\end{theorem}
Theorem $\ref{chap4:thm2}$ says that for an admissible lattice $\tilde{S}(L,\omega,r)$, normalized by a quantity of order $r^{\frac{d-1}{2}}$, converges, in distribution, towards a normal centred distribution. When $L$ is typical, this must not be true: the regularization must be stronger because the orbit $\delta L$ goes repeatedly into the cusp of the space $\mathscr{S}_{d}$, $\textit{id est}$ the zone where $\lVert \delta L \rVert$ is small. \\
In the typical case (typical in the sense of the unique Haar probability measure $\mu_{d}$ over $\mathscr{S}_{d}$), we are going to prove the following result: 
\begin{theorem}
\label{chap4:thm1000}
For every $\epsilon > 0$, for a typical $L \in \mathscr{S}_{d}$, one has that 
\begin{equation}
\label{chap4:eq1032}
\tilde{V}(L,r) = O(r^{2 d - 2 + \epsilon})
\end{equation}
and
\begin{equation}
\label{chap4:eq1033}
\frac{\tilde{V}(L,r)}{r^{ d - 1}} \underset{r \rightarrow \infty}{\rightarrow} \infty \textit{.} 
\end{equation}
In the case where $d=2$, one has also that for every $\epsilon > 0$, for $L$ a typical lattice, one has that
\begin{equation*}
V(L,t) = O(\log(t)^{2+\epsilon})
\end{equation*}
and
\begin{equation*}
\frac{V(L,t)}{\log(t)^{2}} \underset{t \rightarrow \infty}{\rightarrow} \infty \textit{.}
\end{equation*}
Furthermore, in this case, if $P=\text{Rect}(a,b)$, 
$$ \mathbb{E}_{X \in \mathbb{R}^{2}/L}\left( (\frac{\mathcal{R}(t P + X ,L)}{\sqrt{V(L^{\perp},t)}})^{2} \right) \textit{ converges in distribution and in probability towards } \frac{1}{4 \pi^{4} \text{Covol}(L)^{2}} $$
when $t \in [0, T]$ is distributed according to $\frac{1}{T} \rho(\frac{t}{T}) dt$ and when $T \rightarrow \infty$.
\end{theorem}
In particular, in the typical case, the convergence in distribution and in probability of Theorem $\ref{chap4:thm1}$ still holds. Yet, the normalization is larger in this case: "around" $\log(t)$ whereas before it was in $\sqrt{\log(t)}$. \\
Finally, we are also going to tackle another extreme case: the case where $L = \mathbb{Z}^{2}$. $\mathbb{Z}^{2}$ is a unimodular lattice such that there exists $l \in L-\{ 0 \}$ such that $Num(l) = 0$ (for example $(1,0)$). Typically, for $L \in \mathscr{S}_{d}$, $\text{Num}(L)$ is null but there does not exists a non-zero $l \in L$ such that $Num(l) = 0$. \\
In that case, with $\rho$ being a probability density over $[0,1]$, we are going to prove the following theorem: 
\begin{theorem}
\label{chap4:thm3}
For all $x \in \mathbb{R}$, when $t \in [0,T]$ is distributed according to the probability measure $\frac{1}{T} \rho( \frac{t}{T}) dt$ on $[0,T]$ then, when $T \rightarrow \infty$, $\frac{\mathcal{R}(t \text{Rect}(a,a) + (x,x) ,\mathbb{Z}^{2})}{t}$ converges in distribution. Furthermore, the limit distribution $\beta$ has a compact support included in $[-4,4]$ and for every $k \in \mathbb{N}$, one has that 
$$\int_{x \in \mathbb{R}} x^{k} d \beta(x) = a_{k} $$ 
where
$$
a_{k} = \frac{4^{k}(1 + (-1)^{k})(y^{k+1} + (1-y)^{k+1})}{2(k+1)} 
$$
with $y = |t_{2,0} - t_{1,0}|$ where $t_{2,0}$ is the first $t \geqslant 0$ such that $ -t + x \in \mathbb{Z}$ and $t_{1,0}$ is the first $t \geqslant 0$ such that $ t + x \in \mathbb{Z}$.
\end{theorem}
In particular, we see that the normalization in this case is much more important than before and that the error $\mathcal{R}(t \text{Rect}(a,a) + (x,x), \mathbb{Z}^{2})$ in this case is of order $t$. \\
In the next section, we are going to give some heuristic ideas about all these results and then give the plan of the rest of the paper.
\section{Calculation of a Fourier transform, heuristic and plan of the rest of the paper} 
To give some heuristic explanations, we will apply the Poisson formula which states that for a smooth and compact supported function $f: \mathbb{R}^{2} \longrightarrow \mathbb{R}$ and for $L$ a lattice of $\mathbb{R}^{2}$, one has that for every $X \in \mathbb{R}^{2}$ 
\begin{equation}
\label{chap4:eq1100}
\sum_{l \in L} f(l+X) = \frac{1}{\text{Covol}(L)} \sum_{l \in L^{\perp}} \hat{f}(l) e^{-2 i \pi <l,X>} 
\end{equation}
where the Fourier transform $\hat{f}$ is defined by: for every $\xi \in \mathbb{R}^{2}$
\begin{equation}
\label{chap4:eq1101}
\hat{f}(\xi) = \int_{x \in \mathbb{R}^{2}} f(x) e^{2 i \pi <\xi,x>} dx \textit{.}
\end{equation}
In our case, we are interested into the following quantity: 
\begin{equation}
\label{chap4:eq1102}
 N(t P + X,L) = \sum_{l \in L} \mathbf{1}_{t P + X}(l) 
\end{equation}
where $\mathbf{1}_{t P + X}$ is the indicator function of the set $t P + X $ with $P = \text{Rect}(a,b)$. \\
Yet, the function $\mathbf{1}_{t P + X}$ is not smooth enough. But the Poisson formula gives us a good idea of the phenomena that are at play and that is because the Poisson formula applies after having realized a $\textit{smoothing}$ of the studied problem (see Section 4). 
So, we are going to calculate the Fourier transform of $\mathbf{1}_{t P + X}$ in the next subsection. Then we will give the heuristic and announce the plan of the rest of the paper.
\subsection{Calculation of the Fourier transform of $\mathbf{1}_{t P + X}$}
 The main objective of this subsection is to prove the following proposition: 
\begin{prop}
\label{chap4:prop4}
For $P = \text{Rect}(a,b)$, for $l \in \mathbb{R}^{2}$, for $t > 0$, one has: 
$$ \widehat{\mathbf{1}_{t P + X}}(l) = \frac{1}{\pi^{2}}\frac{\sin(2 \pi tl_{1}a)}{ l_{1}} \frac{\sin(2 \pi t l_{2}b)}{ l_{2}}e^{2 i \pi <l,X>} $$
where we convey that $\frac{\sin(0)}{0} = 1$.
\end{prop}
\begin{proof}
Let $l \in \mathbb{R}^{2}$ and $t > 0$. By making the change of variable $x= t y + X$, one has that: 
\begin{equation}
\label{chap4:eq20}
\widehat{\mathbf{1}_{t P + X}}(l) = e^{2 i \pi <l,X>} t^{2} \int_{y \in P} e^{2 i \pi t <l,y>} dy \textit{.}
\end{equation}
Yet, we recall that $P = \text{Rect}(a,b)$ and so one has: 
\begin{equation}
\label{chap4:eq21}
\int_{y \in P} e^{2 i \pi t <l,y>} = \frac{\sin(2 \pi t l_{1} a)}{ \pi t l_{1}}\frac{\sin(2 \pi t l_{2} b)}{ \pi t l_{2}} \textit{.}
\end{equation}
So, with Equation ($\ref{chap4:eq20}$) and Equation ($\ref{chap4:eq21}$), one gets that: 
\begin{equation}
\label{chap4:eq22}
\widehat{\mathbf{1}_{t P + X}}(l) = \frac{e^{2 i \pi <l,X>}}{\pi^{2}} \frac{\sin(2 \pi t l_{1} a)}{l_{1}}\frac{\sin(2 \pi t l_{2} b)}{l_{2}} \textit{.}
\end{equation}
\end{proof}
Now, we are going to give some heuristic explanations about the main results of this paper and the plan of the rest of the paper.
\subsection{Elements of heuristic and plan of the rest of the paper}
Heuristically, the Poisson formula (see Equation $(\ref{chap4:eq1100})$) gives us that: 
\begin{equation}
\label{chap4:eq1054}
 N(t P + X,L) = \sum_{l \in L} \mathbf{1}_{t P + X}(l)  = \frac{1}{\text{Covol}(L)} \sum_{l \in L^{\perp}} \frac{1}{\pi^{2}}\frac{\sin(2 \pi tl_{1}a)}{ l_{1}} \frac{\sin(2 \pi t l_{2}b)}{ l_{2}}e^{2 i \pi <l,X>} 
\end{equation}
with $P = \text{Rect}(a,b)$. \\
Yet, and this link was used a lot in $\cite{trevisan2021}$, in the typical case, the smallest $|l_{1} l_{2}|$ can be seen as a $\lVert \begin{pmatrix} e^{t} & 0 \\ 0 & e^{-t} \end{pmatrix} L \rVert$. This link can also be seen in the following proposition (that will be useful for us later): 
\begin{prop}[\cite{Skriganov}]
\label{chap4:prop1}
$$\text{Num}(L) = d^{- \frac{d}{2}} \inf \{ \lVert \delta L \rVert^{d} \text{ } | \text{ } \delta \in \Delta \}  $$
where $\Delta$ was defined by Equation $\ref{chap4:ajout}$.
\end{prop} 
So it suggests that the term $\frac{\sin(2 \pi tl_{1}a)}{ l_{1}} \frac{\sin(2 \pi t l_{2}b)}{ l_{2}}e^{2 i \pi <l,X>}$ can be as seen as $\frac{\theta_{\delta}(\omega)}{\lVert \delta L \rVert^{d}}$ and so there exists a link between Theorem $\ref{chap4:thm1}$ and Theorem $\ref{chap4:thm2}$. \\ 
Yet, to prove Theorem $\ref{chap4:thm1}$ from Theorem $\ref{chap4:thm2}$, there are two difficulties. The first one is that Theorem $\ref{chap4:thm1}$ and Theorem $\ref{chap4:thm2}$ are about admissible unimodular lattices, which form a negligible set with respect to $\mu_{2}$ and the relation about $|l_{1} l_{2}|$ and $\lVert \begin{pmatrix} e^{t} & 0 \\ 0 & e^{-t} \end{pmatrix} L \rVert$ is not certain in this case. \\
The second one is the fact that the $\sin(2 \pi tl_{1}a)$ and $\sin(2 \pi t l_{2}a)$ do not behave like independent random variable when $T \rightarrow \infty$ and with $t$ being distributed according to the probability measure $\frac{1}{T} \rho(\frac{t}{T}) dt$. Sure we can, like we have done in $\cite{trevisan2021}$ and in $\cite{trevisan2021limit}$, reduce the study to the $l \in L^{\perp}$ that are prime, which is a notion defined by:
\begin{defi}
\label{chap4:def3}
For a lattice $L$, we say that a vector of $L$ is prime if it not a non-trivial integer multiple of another vector of $L$. 
\end{defi}
Yet, even in that case, it is not true that the $\sin(2 \pi tl_{1}a)$ and the $\sin(2 \pi t l_{2}b)$ behave asymptotically like independent random variables. \\
Indeed, let $\alpha \neq \alpha'$ be real irrational numbers with bounded partial quotients in their continued fractions and let us define: 
\begin{equation}
\label{chap4:eq1} \Gamma_{\alpha,\alpha'} = \{ (n + m \alpha, n+ m \alpha' ) \text{ } | \text{ } n,m \in \mathbb{Z}^{2} \} \textit{.}
\end{equation}
Then, we know from $\cite{skriganov1994constructions}$ that $\Gamma_{\alpha,\alpha'}$ is admissible in the sense of Definition $\ref{chap4:def2}$. So, it is also the case of 
\begin{equation}
\label{chap4:eq1049}
 \Gamma_{\alpha,\alpha'}^{\perp} = \frac{1}{\alpha' - \alpha} \{ (n + m \alpha, n+ m \alpha' ) \text{ } | \text{ } n,m \in \mathbb{Z}^{2} \} \textit{.}
\end{equation}
Then, if we consider, for $k \geqslant 2$, 
\begin{equation}
\label{chap4:eq1050}
v_{1}(k) = \frac{1}{\alpha' - \alpha} (k + (k+1) \alpha, k + (k+1) \alpha' ) \textit{, }
\end{equation}
\begin{equation}
\label{chap4:eq1051}
v_{2}(k) = \frac{1}{\alpha' - \alpha} (k+1 + (k+2) \alpha, k+1 + (k+2) \alpha' ) \textit{, }
\end{equation}
\begin{equation}
\label{chap4:eq1052}
v_{3}(k) = \frac{1}{\alpha' - \alpha} (k+2 + (k+3) \alpha, k+2 + (k+3) \alpha' ) \textit{, }
\end{equation}
we have that: $v_{1}(k),v_{2}(k)$ and $v_{3}(k)$ are prime vectors of $ \Gamma_{\alpha,\alpha'}^{\perp}$ and
\begin{equation}
\label{chap4:eq1053}
-v_{1}(k) + 2 v_{2}(k) = v_{3}(k)
\end{equation}
which prevents from getting the wanting asymptotic independence. \\
That is why, to get rid of this problem, we consider $\mathbb{E}_{X \in \mathbb{R}^{2}/L}\left( (\frac{\mathcal{R}(t P + X ,L)}{\sqrt{V(L^{\perp},t)}})^{2} \right)$ instead of considering $\frac{\mathcal{R}(t P + X ,L)}{\sqrt{V(L^{\perp},t)}}$. \\
By doing so, the Parseval formula gives us that, if we want to prove Theorem $\ref{chap4:thm1}$, heuristically we have to study the convergence in distribution of 
\begin{equation}
\label{chap4:eq1055}
\frac{G(L^{\perp},t)}{V(L^{\perp},t)}  = \frac{2}{\pi^{4} \text{Covol}(L)^{2} V(L^{\perp},t) } \sum_{l \in J_{2}(L^{\perp},t)} \frac{1}{(l_{1}l_{2})^{2}} \sum_{k=1}^{\infty} \frac{(\sin(2 \pi k t l_{1}) \sin(2 \pi k t l_{2}))^{2} }{k^{4}}
\end{equation}
where 
\begin{equation}
\label{chap4:eq1200}
J_{2}(L^{\perp},t) = \{ l \in L^{\perp} \text{ } | \text{ } 0 < \lVert l \rVert \leqslant t \text{ , } l \text{ prime} \text{ and } l_{1} > 0 \} 
\end{equation}
and where $t \in [0, T]$ is being distributed according to the probability measure $\frac{1}{T} \rho(\frac{\cdot}{T}) dt$ and where $T \rightarrow \infty$. We have implicitly used Equation $(\ref{chap4:eq1054})$ and centred the sum of the right-hand side on the prime vectors whose norm are smaller than $t$. We have to cut the sum because of a problem of convergence. \\
The final ideas that we use to prove Theorem $\ref{chap4:thm1}$ are, first, the idea to use the well-known formula $\sin^{2}(\cdot) = \frac{1-\cos(2 \cdot)}{2}$. Then there are two different types of quantities that must be dealt with. \\
The first one is of the type:
\begin{equation}
\label{chap4:eq1056}
\frac{1}{V(L^{\perp},t)} \sum_{\substack{ l \in L^{\perp}-\{0 \} \\ \lVert l \rVert \leqslant t }} \frac{1}{ l_{1}^{2} l_{2}^{2} }
\end{equation}
and we show quickly that such a quantity converges almost surely when $T \rightarrow \infty$. \\
The other type of quantities have a term of the form $\cos(t f(l))$ or a product of two $\cos(t f(l))$ (with $f$ being a function of $l$) in the numerator of the terms. In that case, we show that the moment of order $2$ of the quantities of this type converge to $0$ when $T \rightarrow \infty$. \\
We use the moment of order $2$ because it is quite convenient since we are dealing with numerators that have a term of the form $\cos(t f(l))$. Yet, this last part is a calculatory one. Furthermore, we have to underline the fact that these calculations still work for the typical $L$ considered in Theorem $\ref{chap4:thm1000}$, not only for admissible lattice $L$. They are, in that sense, intrinsic. %à vérifier \\
The estimate of $V(L,t)$ contained in Theorem $\ref{chap4:thm1}$ is basically derived from the estimate of the integral $$ \int_{\substack{ l \in \mathbb{R}^{2} \\ A \leqslant \lVert l \rVert \leqslant t \\ |\text{Num}(l)| \geqslant C}} \frac{1}{(l_{1}l_{2})^{2}} $$ where $A > 0$, $C > 0$ and $t \rightarrow \infty$. It concludes the heuristic explanations that we wanted to give about Theorem $\ref{chap4:thm1}$.\\
\\
The proof of Theorem $\ref{chap4:thm2}$ is quicker and it is basically an application of the central limit theorem with error term (see Theorem $\ref{chap4:thm4}$). \\
We already have said a few words about the last part of Theorem $\ref{chap4:thm1000}$. We will now give some explanations about the estimates $V$ and $\tilde{V}$ presented in this theorem. \\
The estimates of $\tilde{V}$ are applications of the ergodic theorem whereas the estimates of $V$ is deduced from an upper estimate of 
$$\int_{\substack{ l \in \mathbb{R}^{2} \\ A \leqslant \lVert l \rVert \leqslant t  \\ |\text{Num}(l)| \geqslant C | \log(\lVert l \rVert) |^{-1-\alpha}}} \frac{1}{l_{1}^{2}  l_{2}^{2}} dl_{1} dl_{2} $$
where $C > 0$, $A > 0$ and $\alpha > 0$ with $t \rightarrow \infty$ and from results of $\cite{sprindzhuk1979metric}$ and of $\cite{Skriganov}$ (see Theorem $\ref{chap4:thm12}$). \\
Finally, concerning the Theorem $\ref{chap4:thm3}$, the appropriate normalization of the error term $\mathcal{R}$ is $t$ because, to within a multiplicative factor, it is the perimeter of $t \text{Rect}(a,a) + X$. Indeed, in this case, we can easily compute $\mathcal{R}$. The wanted convergence of Theorem $\ref{chap4:thm3}$ is then deduced from it and from the application of the method of moments. \\
\\
$\textbf{Plan of the paper.}$ In Section 3, we prove Theorem $\ref{chap4:thm2}$. \\ In Section 4, we prove Theorem $\ref{chap4:thm1}$. The first subsection is dedicated to obtain the upper and lower estimate on $V(L,t)$ in the case where $L$ is admissible. \\
Let us set: $$ S(L,X,t) = \frac{2}{\pi^{2} \textit{covol}(L)} \sum_{l \in J_{2}(L^{\perp},t)} \frac{1}{l_{1}l_{2}} \sum_{k=1}^{\infty} \frac{\sin(2 \pi k t l_{1}) \sin(2 \pi t l_{2}) \cos(2 \pi k <l,X> ) }{k^{2}} $$ 
where we recall that $J_{2}(L^{\perp},t) = \{ l \in L^{\perp} \text{ } | \text{ } 0 < \lVert l \rVert \leqslant t \text{ , } l \text{ prime} \text{ and } l_{1} > 0 \} $. \\
 The second subsection is dedicated to show that we can reduce the study of the convergence in distribution of $\mathbb{E}_{X \in \mathbb{R}^{2}/L}\left( (\frac{\mathcal{R}(t P + X ,L)}{\sqrt{V(L^{\perp},t)}})^{2} \right)$ to the study of the convergence in distribution of $\mathbb{E}_{X \in \mathbb{R}^{2}/L}\left( \frac{S(L,X,t)}{\sqrt{V(L^{\perp},t)}}^{2} \right)$ when $t$ is being distributed according to $\frac{1}{T} \rho(\frac{t}{T}) dt$ and when $T \rightarrow \infty$ (see Proposition $\ref{chap4:prop3}$). The third subsection is dedicated to show that $ \mathbb{E}_{X \in \mathbb{R}^{2}/L}\left( \frac{S(L,X,t)}{\sqrt{V(L^{\perp},t)}}^{2} \right)$ (see Proposition $\ref{chap4:prop8}$) converges in distribution when $t$ is being distributed according to $\frac{1}{T} \rho(\frac{t}{T}) dt$ and when $T \rightarrow \infty$. The fourth subsection concludes the proof of Theorem $\ref{chap4:thm1}$.  \\
In Section 5, we give the proof of Theorem $\ref{chap4:thm1000}$. The first subsection is dedicated to the estimates of $\tilde{V}(L,r)$. The second subsection is dedicated to the estimates of $V(L,t)$. The third subsection concludes the proof of Theorem $\ref{chap4:thm1000}$, using, in particular, the third subsection of Section 4. \\
In Section 6, we give the proof of Theorem $\ref{chap4:thm3}$. In the first subsection, we give a simple expression of $\frac{\mathcal{R}(t \text{Rect}(1,1) + (x,x) ,\mathbb{Z}^{2})}{t}$ (and before we see that the $a$, in the statement of Theorem $\ref{chap4:thm3}$, can be chosen equal to $1$). In the second subsection, we can reduce the study of the convergence in distribution of $ \frac{\mathcal{R}(t \text{Rect}(1,1) + (x,x) ,\mathbb{Z}^{2})}{t}$, when $T \rightarrow \infty$, to a simpler quantity. In the third subsection, we are going to apply the method of moments to get the convergence of this simpler quantity. In the fourth subsection, we conclude the proof of Theorem $\ref{chap4:thm3}$. 
\section{Proof of Theorem 14}
To prove Theorem $\ref{chap4:thm2}$, we need to recall two important theorems. The first one is the central limit theorem with error term (see for example $\cite{beck2010randomness}$): 
\begin{theorem}
\label{chap4:thm4}
Let $(Z_{i})_{1 \leqslant i \leqslant n}$ be a sequence of real random variables independent such that for all $1 \leqslant i \leqslant n $, $\mathbb{E}(Z_{i}) = 0$ and $Z_{i}$ admits a moment of order $3$. Let us call: 
$$ T = \sum_{i=1}^{n} \mathbb{E}(|Z_{i}|^{3}) \textit{ and } V = \sum_{i=1}^{n} \mathbb{E}(Z_{i}^{2}) \textit{.} $$
Then for every real $\lambda$, one has 
$$| \mathbb{P}( \frac{\sum_{i=1}^{n}Z_{i}}{\sqrt{V}} \geqslant \lambda ) - \frac{1}{\sqrt{2 \pi}} \int_{\lambda}^{\infty} e^{- \frac{u^{2}}{2}} du| < \frac{40 T}{V^{\frac{3}{2}}} \textit{.}$$
\end{theorem}
The second one is a theorem due to Minkowski: 
\begin{theorem}
\label{chap4:thm5}
There exists $K > 0 $  such that for all $L \in \mathscr{S}_{d}$, 
$$\lVert L \rVert  \leqslant K \textit{.}$$
\end{theorem}
We can now prove Theorem $\ref{chap4:thm2}$.
\begin{proof}[Proof of Theorem $\ref{chap4:thm2}$]
Let $L$ be an admissible lattice. So, according to Proposition $\ref{chap4:prop1}$ and Equation $\ref{chap4:ajout}$, there exists $C > 0$ such that for every $\delta \in \Delta$, one has 
\begin{equation}
\label{chap4:eq2}
\lVert \delta L \rVert \geqslant C \textit{.}
\end{equation}
Theorem $\ref{chap4:thm5}$ gives us that: 
\begin{equation}
\label{chap4:eq3}
\lVert \delta L \rVert \leqslant K \textit{.}
\end{equation}
Because of Equation ($\ref{chap4:eq2}$) and Equation ($\ref{chap4:eq3}$), one has  
\begin{equation}
\label{chap4:eq4}
| \Delta_{r} | K^{-2d} \leqslant \tilde{V}(L,r) \leqslant | \Delta_{r} | C^{-2d} \textit{.}
\end{equation}
Yet, the cardinal number of $\Delta_{r}$ satisfies that 
\begin{equation} 
\label{chap4:eq1300}
|\Delta_{r}| = n_{r}r^{d-1} + o_{r}(r^{d-1}) \textit{.}
\end{equation}
So we get the first wanted result. \\
According to Theorem $\ref{chap4:thm5}$, one has: 
\begin{equation}
\label{chap4:eq5}
| \mathbb{P}( \frac{\tilde{S}(L,\omega,r)}{\sqrt{ \tilde{V}(L,r)}} \geqslant \lambda ) - \frac{1}{\sqrt{2 \pi}} \int_{\lambda}^{\infty} e^{- \frac{u^{2}}{2}} du| < \frac{40 T_{1}(r)}{V_{1}(r)^{\frac{3}{2}}}
\end{equation}
where 
\begin{equation}
\label{chap4:eq6}
V_{1}(r) = \sum_{\delta \in \Delta_{r}} \mathbb{E}((\frac{\theta_{\delta}}{\lVert \delta L  \rVert^{d}})^{2})
\end{equation}
and
\begin{equation}
\label{chap4:eq7}
T_{1}(r) = \sum_{\delta \in \Delta_{r}} \mathbb{E}(|\frac{\theta_{\delta}}{\lVert \delta L  \rVert^{d}}|^{3}) \textit{.}
\end{equation}
Yet, one has: 
\begin{equation}
\label{chap4:eq8}
V_{1}(r) = \Theta(r^{d-1})
\end{equation}
and
\begin{equation}
\label{chap4:eq9}
T_{1}(r) = \Theta(r^{d-1})
\end{equation}
because the $\theta_{\delta}$ are independent, symmetrical, identically distributed and their common distribution admit a moment of order 3 and because of Equation ($\ref{chap4:eq2}$) and of Equation ($\ref{chap4:eq3}$). \\
Thus, one gets that: 
\begin{equation}
\label{chap4:eq10}
| \mathbb{P}( \frac{\tilde{S}(L,\omega,r)}{\sqrt{ \tilde{V}(L,r)}} \geqslant \lambda ) - \frac{1}{\sqrt{2 \pi}} \int_{\lambda}^{\infty} e^{- \frac{u^{2}}{2}} du| = O(\frac{1}{r^{\frac{d-1}{2}}}) \textit{.}
\end{equation}

From this last Equation and by making $r \rightarrow \infty$, one gets the wanted result.
\end{proof}
%premièrement le théorème de convergence en loi vers loi normale 
%deuxièmement Minkowksi 
\section{Proof of Theorem 13}
\subsection{Estimates of $V(L,t)$ when $L$ is admissible}
%Dire quelques mots d'explication, proposition principale ...
In this subsection, we are going to prove the following proposition, which is the first assertion of Theorem $\ref{chap4:thm1}$ 
\begin{prop}
\label{chap4:prop2}
Let $L$ be an admissible lattice of $\mathbb{R}^{2}$. One has that there exists $C > 0$ (big enough) such that for all $t$ large enough: 
$$\frac{1}{C} \log(t) \leqslant V(L,t) \leqslant C \log(t) \textit{.} $$
\end{prop}
The main idea is that $V(L,t)$ can be compared to
$$ \int_{\substack{ l \in \mathbb{R}^{2} \\ A \leqslant \lVert l \rVert \leqslant t+B \\ |\text{Num}(l)| \geqslant C}} \frac{1}{l_{1}^{2}  l_{2}^{2}} dl_{1} dl_{2} $$
when $L$ is admissible and when $C,A$ and $B$ are positive constants that are well-chosen. And this last integral behaves like $\log(t)$ to within one multiplicative constant. \\
This last fact is the object of the following lemma.
%\begin{lemma}
%\label{chap4:lemme1}
%Let $f: \mathbb{R}^{d} \longrightarrow \mathbb{R}_{+}-\{  0 \}$ be a measurable function that is decreasing in each of its variables. Then, there exists $C > 0 $ such that: 
%$$\frac{1}{C} \int_{\mathbb{R}^{d}} f(x) dx \leqslant \sum_{l \in L} f(l) \leqslant C \int_{\mathbb{R}^{d}} f(x) dx \textit{.}$$
%\end{lemma}
In this lemma, we use the notation $$f(t) \sim_{t \rightarrow \infty} g(t) $$ to express the facts that, for $t$ large enough, $g(t) \neq 0$ and that $\frac{f(t)}{g(t)} \rightarrow 1$ when $t \rightarrow \infty$. \\
%The second one is an estimating one: 
\begin{lemma}
\label{chap4:lemme2}
For all $C > 0$, for all $A > 0$, for all $t$ large enough, one has: 
$$\int_{\substack{ l \in \mathbb{R}^{2} \\ A \leqslant \lVert l \rVert \leqslant t \\ |\text{Num}(l)| \geqslant C}} \frac{1}{l_{1}^{2}  l_{2}^{2}} dl_{1} dl_{2} \sim_{t \rightarrow \infty} \frac{8 \log(t)}{ C }   \textit{.}$$
\end{lemma}
\begin{proof}
Let us set: 
\begin{equation}
J(t) = \int_{\substack{ l \in \mathbb{R}^{2} \\ A \leqslant \lVert l \rVert \leqslant t \\ |\text{Num}(l)| \geqslant C}} \frac{1}{l_{1}^{2}  l_{2}^{2}} dl_{1} dl_{2}  
\end{equation}
and let us remark by the way that it is enough to prove the result for $A$ large enough. \\
By passing into polar coordinates $(r,\theta)$ and by using the symmetries, one has: 
\begin{equation}
\label{chap4:eq11}
J(t) = 8 \int_{\substack{A \leqslant r \leqslant t \\ 0 \leqslant \theta \leqslant \frac{\pi}{4} \\ \sin(2 \theta) \geqslant \frac{2 C }{r^{2}} }} \frac{1}{r^{3}} \frac{4}{\sin(2 \theta)^{2}} dr d\theta \textit{.}
\end{equation} %16
By making the changes of variable $\theta' = 2 \theta$ and, then, $u = \tan(\theta')$, one gets from Equation $(\ref{chap4:eq11})$ and by taking $A$ large enough: 
\begin{equation}
\label{chap4:eq12}
J(t) = 16 \int_{A \leqslant r \leqslant t} \frac{1}{r^{3}} (-1 + \frac{r^{2}}{2 C} \sqrt{1 - (\frac{2 C}{r^{2}})^{2}}) dr \textit{.} 
\end{equation}
From Equation ($\ref{chap4:eq12}$) and because of the fact that $ \frac{1}{r^{3}} (-1 + \frac{r^{2}}{2 C} \sqrt{1 - (\frac{2 C}{r^{2}})^{2}}) =  \frac{1}{2 r  C } + O(\frac{1}{r^{2}})$ (when $r \rightarrow \infty$), we finally get that: 
\begin{equation}
\label{chap4:eq13}
J(t) \sim_{t \rightarrow \infty} \frac{8}{ C } \log(t) \textit{.}
\end{equation}
So we get the wanted result.
\end{proof}
We can now prove Proposition $\ref{chap4:prop2}$: 
\begin{proof}[Proof of Proposition $\ref{chap4:prop2}$]
%There exists $D > 0$ such that for all $t$ big enough: 
%\begin{equation}
%\label{chap4:eq18}
%\frac{1}{D} \sum_{ \lVert l \rVert \leqslant t} \frac{1}{l_{1}^{2} l_{2}^{2} + K} \leqslant H(t) \leqslant D \sum_{ \lVert l \rVert \leqslant t} \frac{1}{l_{1}^{2} l_{2}^{2} + K} 
%\end{equation}
%where $K > 0$ is such that for all $l \in L-\{ 0 \}$, $|l_{1}l_{2}| \geqslant K$. K exists because $L$ is supposed to be admissible. \\
%Lemme $\ref{chap4:lemme1}$ gives us that: 
%\begin{equation}
%\label{chap4:eq15}
%\sum_{ \lVert l \rVert \leqslant t} \frac{1}{l_{1}^{2} l_{2}^{2} + K} = \Theta(\int_{ \lVert l \rVert \leqslant t } \frac{1}{K + l_{1}^{2} l_{2}^{2}} dl_{1} dl_{2})
%\end{equation} 
%and furthermore, one has: 
%\begin{equation}
%\label{chap4:eq16}
%%\int_{ \lVert l \rVert \leqslant t } \frac{1}{K + l_{1}^{2} l_{2}^{2}} dl_{1} dl_{2} = \Theta(\int_{ \lVert l \rVert \leqslant t } \frac{1}{l_{1}^{2} l_{2}^{2}} dl_{1} dl_{2}) \textit{.}
%\end{equation} 
First, as $L$ is admissible, there exists $C >0$ such that for all $l \in L-\{ 0 \}$, $$|l_{1}l_{2}| \geqslant C \textit{.} $$
So, according to Definition $\ref{chap4:def10}$ and by integration by parts, there exists $A > 0$, $B > 0$ and $D > 0$ such that: 
\begin{equation}
\label{chap4:eq1057}
\frac{1}{D} \int_{\substack{ l \in \mathbb{R}^{2} \\ A \leqslant \lVert l \rVert \leqslant t+B \\ |\text{Num}(l)| \geqslant C}} \frac{1}{l_{1}^{2}  l_{2}^{2}} dl_{1} dl_{2}  \leqslant V(L,t) \leqslant D \int_{A \leqslant \lVert l \rVert \leqslant t+B } \min(\frac{1}{C^{2}}, \frac{1}{(l_{1}l_{2})^{2}}) dl_{1}dl_{2} \textit{.}
\end{equation}
Lemma $\ref{chap4:lemme2}$ gives us in particular that: 
\begin{equation}
\label{chap4:eq17}
\int_{\substack{ l \in \mathbb{R}^{2} \\ A \leqslant \lVert l \rVert \leqslant t+B \\ |\text{Num}(l)| \geqslant C}} \frac{1}{l_{1}^{2}  l_{2}^{2}} dl_{1} dl_{2} = \Theta(\log(t)) \textit{.}
\end{equation} 
and
\begin{equation}
\label{chap4:eq1058}
\int_{A \leqslant \lVert l \rVert \leqslant t+B } \min(\frac{1}{C^{2}}, \frac{1}{(l_{1}l_{2})^{2}}) dl_{1}dl_{2} = \Theta(\log(t)) + \frac{1}{C^{2}} \int_{\substack{ A \leqslant \lVert l \rVert \leqslant t+B \\ |Num(l)| \leqslant C}} dl_{1} dl_{2} \textit{.}
\end{equation} 
Yet, one has that: 
\begin{equation}
\label{chap4:eq1059}
\int_{\substack{ A \leqslant \lVert l \rVert \leqslant t+B \\ |Num(l)| \leqslant C}} dl_{1} dl_{2} = \Theta(\log(t)) \textit{.}
\end{equation}
Equation ($\ref{chap4:eq1057}$), Equation ($\ref{chap4:eq17}$), Equation ($\ref{chap4:eq1058}$) and Equation ($\ref{chap4:eq1059}$) give us then the wanted result.
\end{proof}

$\textit{Remark.}$ We think that in dimension $d \geqslant 3$, we must have that for all $C > 0$, for all $A > 0$, one has: 
\begin{equation}
\label{chap4:eq14}
\int_{\substack{ l \in \mathbb{R}^{d} \\ A \leqslant \lVert l \rVert \leqslant t \\ |\text{Num}(l)| \geqslant C}} \frac{1}{l_{1}^{2} \cdots  l_{d}^{2}} dl_{1} \cdots dl_{d} = \Theta(\log(t)^{d-1}) \textit{.}
\end{equation}
In fact, the upper part of Equation ($\ref{chap4:eq14}$) can be proven like that: 
\begin{proof}
By symmetry, we are brought back to the case where $l_{i} > 0$ for all $i$.\\
We set $$ \phi: (l_{1}, \cdots , l_{d}) \longmapsto (l_{1}, l_{1}l_{2}, \cdots, l_{1} \cdots l_{d}) \textit{.}$$ Then $\phi$ is a $C^{\infty}$-diffeomorphism from $(\mathbb{R}_{+}-\{ 0 \})^{d}$ to itself and the jacobian matrix evaluated on $(l_{1}, \cdots, l_{d})$, denoted by $\text{Jac}(l_{1}, \cdots , l_{d})$, satisfies that: 
\begin{equation}
\label{chap4:eq153}
\text{Jac}(l_{1}, \cdots , l_{d}) = \prod_{i=1}^{d-1} \phi_{i}(l_{1},\cdots,l_{d}) \textit{.}
\end{equation}
Furthermore, because $l$ belongs to the domain of the integral, we see that for all $i \in \{1, \cdots, d-1 \}$, $$ b_{i} = \frac{1}{C t^{(1+\epsilon_{1})(d-i)}} \leqslant \phi_{i}(l) \leqslant t^{(1+\epsilon_{1})i}=h_{i} \textit{ and}$$ $$\phi_{d}(l) \geqslant C =b_{d} \textit{.} $$
Hence, from ($\ref{chap4:eq153}$), we get that: 
\begin{equation}
\label{chap4:eq154}
\int_{\substack{\forall 1 \leqslant i \leqslant d \textit{, } l_{i} > 0 \ \\ A \leqslant \lVert l \rVert_{2} \leqslant t \\ \text{Num}(l) \geqslant C}} \frac{1}{l_{1}^{2} \cdots l_{d}^{2}} dl_{1} \cdots dl_{d} \leqslant \int_{ \substack{ \forall 1 \leqslant i \leqslant d-1 \textit{, } b_{i} \leqslant u_{i} \leqslant h_{i} \\ b_{d} \leqslant u_{d} }} \frac{1}{u_{d}^{2} \prod_{i=1}^{d-1} u_{i}} du \textit{.} 
\end{equation}
The right hand side can be easily calculated and gives the wanted result.
\end{proof}
To prove the lower part of the Equation $(\ref{chap4:eq14})$, we think that one way is to use hyperspherical coordinates. \\
\\
We are now going to tackle the proof of the second part of Theorem $\ref{chap4:thm1}$. It should be noted that, by density, it is enough to treat the case where the support of $\rho$ is included in $[\alpha,1]$ with $0 < \alpha < 1$.
\subsection{Smoothing and reduction to the study of a Fourier series}
In this subsection, we are going to prove the following proposition: 
\begin{prop}
\label{chap4:prop3}
Let $L$ be an admissible lattice of $\mathbb{R}^{2}$. \\
Then, for all $T$ large enough,  we have that: 
$$ \mathbb{E}_{X \in \mathbb{R}^{2}/L} \left( (\frac{\mathcal{R}(  t \mathbb{D}^{2}+X,L)}{\sqrt{V(L^{\perp},t)})}-\frac{S(L^{\perp},X,t)}{\sqrt{V(L^{\perp},t)}})^{2}  \right) \rightarrow 0$$
where the convergence towards $0$ is uniform in $\alpha T \leqslant t \leqslant T$ when $T \rightarrow \infty$, where
\begin{equation}
\label{chap4:eq19}
S(L,X,t) = \frac{2}{\pi^{2} \textit{covol}(L)} \sum_{l \in J_{2}(L,t)} \frac{1}{l_{1}l_{2}} \sum_{k=1}^{\infty} \frac{\sin(2 \pi k t l_{1}) \sin(2 \pi k t l_{2}) \cos(2 \pi k <l,X> ) }{k^{2}} 
\end{equation}
and where 
\begin{equation}
\label{chap4:eq80}
J_{2}(L,t) = \{ l \in L \text{ } | \text{ } 0 < \lVert l \rVert \leqslant t \text{, } l \text{ prime} \text{ and } l_{1} > 0 \} \textit{.}
\end{equation}
\end{prop}
Proposition $\ref{chap4:prop3}$ basically tells us that the asymptotical study of $ \mathbb{E}_{X \in \mathbb{R}^{2}/L}\left( (\frac{\mathcal{R}(t P + X ,L)}{\sqrt{V(L^{\perp},t)}})^{2} \right)$, when $T \rightarrow \infty$, can be reduced to the study of $ \mathbb{E}_{X \in \mathbb{R}^{2}/L}\left( \frac{S(L^{\perp},X,t,T)}{\sqrt{V(L^{\perp},t)}} \right)$ (it is in fact stronger). It is suggested by the Poisson formula. Yet, we can not use it directly to prove this fact because the indicator function $\mathbf{1}_{t P +X}$, with $$P = \textit{Rect}(a,b)$$ is not regular enough. \\
So, the first thing we are going to do to prove Proposition $\ref{chap4:prop3}$ is to $\textit{smooth}$ the studied problem. It is the object of the next subsubsection. The study of $\frac{\mathcal{R}(t P + X ,L)}{\sqrt{V(L^{\perp},t)}}$ is going to be reduced to the study of two Fourier series. Then, we will do some calculations to simplify this study and finally the study of $\frac{\mathcal{R}(t P + X ,L)}{\sqrt{V(L^{\perp},t)}}$ is going to be reduced to the study of $\frac{S(L^{\perp},X,t)}{\sqrt{V(L^{\perp},t)}}$.
%\subsection{Calculation of the Fourier transform of the characteristic function of $ \text{Rect}(a,b)$}
%The main object is to prove the following proposition: 
%\begin{prop}
%\label{chap4:prop4}
%For $P = \text{Rect}(a,b)$, for $l \in \mathbb{R}^{2}$, for $t > 0$, one has: 
%$$ \widetilde{\mathbf{1}_{t P + X}}(l) = \frac{1}{\pi^{2}}\frac{\sin(2 \pi tl_{1}a)}{ l_{1}} \frac{\sin(2 \pi t l_{2}a)}{ l_{2}}e^{2 i \pi <l,X>} $$
%where we convey that $\frac{\sin(0)}{0} = 1$.
%\end{prop}
%\begin{proof}
%Let $l \in \mathbb{R}^{2}$ and $t > 0$. By making the change of variable $x= t y + X$, one has that: 
%\begin{equation}
%\label{chap4:eq20}
%\widetilde{\mathbf{1}_{t P + X}}(l) = e^{2 i \pi <l,X>} t^{2} \int_{y \in P} e^{2 i \pi t <l,y>} dy \textit{.}
%\end{equation}
%Yet, we recall that $P = \text{Rect}(a,b)$ and so one has: 
%\begin{equation}
%\label{chap4:eq21}
%\int_{y \in P} e^{2 i \pi t <l,y>} = \frac{\sin(2 \pi t l_{1} a)}{ \pi t l_{1}}\frac{\sin(2 \pi t l_{2} b)}{ \pi t l_{2}} \textit{.}
%\end{equation}
%So, with equations ($\ref{chap4:eq20}$) and ($\ref{chap4:eq21}$), one gets that: 
%\begin{equation}
%\label{chap4:eq22}
%\widetilde{\mathbf{1}_{t P + X}}(l) = \frac{e^{2 i \pi <l,X>}}{\pi^{2}} \frac{\sin(2 \pi t l_{1} a)}{l_{1}}\frac{\sin(2 \pi t l_{2} b)}{l_{2}} \textit{.}
%\end{equation}
%\end{proof}
%With this last formula, we see, again, that we need to $\textit{regularize}$ the problem if we want to use the Poisson formula. The next subsection is dedicated to it. 
\subsubsection{Reduction to the study of two Fourier series}
%suivre Skriganov ici
This subsubsection is dedicated to the $\textit{smoothing}$ of the problem. We are going to do it like in $\cite{Skriganov}$.\\
Let us take $\omega$ a $C^{\infty}$-function, of compact support included in $B_{f}(0,1)$, such that $$\int_{\mathbb{R}^{2}} \omega(x) dx = 1 \textit{, }$$ $\omega(x) \geqslant 0$ and such that $\omega$ is spherically symmetric (so will be its Fourier transform). \\
Let us recall the following definition: 
\begin{defi}
\label{chap4:def4}
Let $\mathbf{O} \subset \mathbb{R}^{2}$ be a connected compact set of $\mathbb{R}^{2}$ with a boundary that is regularly piecewise. Let $1 > \tau > 0$. A couple of compact region $(\mathbf{O}_{\tau}^{+},\mathbf{O}_{\tau}^{-})$ is called a $\tau$-co-approximation if $ \mathbf{O}_{\tau}^{-} \subset \mathbf{O} \subset \mathbf{O}_{\tau}^{+} $ and if the points at the frontiers $\partial \mathbf{O}_{\tau}^{\pm}$ are at least distant from $\tau$ of the frontier $\partial \mathbf{O}$. 
\end{defi}
For example, we can think that $\mathbf{O} = P$ and, in that case, a $\tau$-co-approximation is given by $(1 \pm \beta \tau) P$, with $\beta > 0$ well-chosen. \\
Let us define: 
\begin{equation}
\label{chap4:eq23}
\mathfrak{R}_{\tau}^{\pm}(\mathbf{O},X) = \frac{1}{covol(L)} \sum_{l \in L^{\perp}-\{ 0\}} \hat{\mathbf{1}}_{\mathbf{O}_{\tau}^{\pm}}(l) \hat{\omega}(\tau l) e^{ 2 i \pi <l,X> } \textit{.}
\end{equation}
For $\mathbf{O} = t P = t \text{Rect}(a,b)$, there exists $\beta > 0$ (independent from $t$) such that $\mathbf{O}_{\tau}^{\pm} = (t \pm \beta \tau ) P$ is a $\tau$-co-approximation. \\
Let us take: \begin{equation}
\label{chap4:eq28}
\tau = \tau(t) = \frac{\log(t)^{\gamma}}{t} 
\end{equation}
where $ \frac{1}{2} > \gamma > 0$. \\
The main objective of this subsubsection is to prove the following proposition: 
\begin{prop}
\label{chap4:prop5}
For this choice of $\tau$ and of $\tau$-co-approximation $(t \pm \beta \tau )P$, there exists $B > 0$ (independent from X) such that for all $T$ large enough 
$$
 \frac{\mathfrak{R}_{\tau}^{-}(tP,X)}{\sqrt{V(L^{\perp},t)}} - B \log(T)^{\gamma - \frac{1}{2}} \leqslant \frac{\mathcal{R}(tP+X,L)}{\sqrt{V(L^{\perp},t)}} \leqslant \frac{\mathfrak{R}_{\tau}^{+}(tP,X)}{\sqrt{V(L^{\perp},t)}} + B\log(T)^{\gamma - \frac{1}{2}} \textit{.}
$$
\end{prop}
Basically, it says that one can brought back the study of $\frac{\mathcal{R}(tP+X,L)}{\sqrt{V(L^{\perp},t)}}$ to the study of the two quantities $\frac{\mathfrak{R}_{\tau}^{\pm}(tP,X)}{\sqrt{V(L^{\perp},t)}}$, which represent the $\textit{smoothed}$ Fourier series of $\frac{\mathcal{R}(tP+X,L)}{\sqrt{V(L^{\perp},t)}}$. \\
To prove this proposition, the following lemma is going to be useful: 
\begin{lemma}
\label{chap4:lemme3}
Let us take $0 < \tau < 1$. Then, one has: 
$$\mathfrak{R}_{\tau}^{-}(\mathbf{O},X) + \frac{\text{Area}(\mathbf{O}_{\tau}^{-}) - \text{Area}(\mathbf{O})}{covol(L)} \leqslant \mathcal{R}(\mathbf{O}+X,L) \leqslant \mathfrak{R}_{\tau}^{+}(\mathbf{O},X) + \frac{\text{Area}(\mathbf{O}_{\tau}^{+}) - \text{Area}(\mathbf{O})}{covol(L)} \textit{.}$$
\end{lemma}
\begin{proof}
Let us set $\omega_{\tau}(x) = \tau^{-2} \omega(\tau^{-1} x)$ for all $x \in \mathbb{R}^{2}$. Then, one has: 
\begin{equation}
\label{chap4:eq24}
\hat{\omega}_{\tau}(\xi) = \hat{\omega}(\tau \xi) 
\end{equation}
for all $\xi \in \mathbb{R}^{2}$. \\
Let us consider the following convolution products: 
\begin{equation}
\label{chap4:eq25}
h^{\pm}(x) = (\omega_{\tau} * \mathbf{1}_{\mathbf{O}_{\tau}^{\pm}})(x) = \int_{\mathbb{R}^{2}} \omega_{\tau}(x-y) \mathbf{1}_{\mathbf{O}_{\tau}^{\pm}}(y) dy \textit{.} 
\end{equation}
The functions $h^{\pm}$ are $C^{\infty}$ over $\mathbb{R}^{2}$ and have a compact support. \\
By using Definition $\ref{chap4:def4}$, one has, for all $x \in \mathbb{R}^{2}$: 
\begin{equation}
\label{chap4:eq26}
h^{-}(x) \leqslant \mathbf{1}_{\mathbf{O}}(x) \leqslant h^{+}(x) \textit{.} 
\end{equation}
By replacing $x \in \mathbb{R}^{2}$ by $l-x$ with $l \in L$ and by summing over $l \in L$, one gets with Equation ($\ref{chap4:eq26}$): 
\begin{equation}
\label{chap4:eq27}
\sum_{l \in L} h^{-}(l-X) \leqslant N( \mathbf{O} + X, L) \leqslant \sum_{l \in L} h^{+}(l-X) \textit{.}
\end{equation}
By using the Poisson formula and the fact that Fourier transform a convolution product into an usual product, one gets the wanted result thanks to Equation ($\ref{chap4:eq24}$). 
\end{proof}
We can now tackle the proof of Proposition $\ref{chap4:prop5}$: 
\begin{proof}[Proof of Proposition $\ref{chap4:prop5}$]
According to Lemma $\ref{chap4:lemme3}$, one has that, for every $t>0$ and $X$,
\begin{equation}
\label{chap4:eq29}
\mathfrak{R}_{\tau}^{-}(tP,X) + \frac{\text{Area}((t - \beta \tau ) P) - \text{Area}(t P)}{covol(L)} \leqslant \mathcal{R}(tP+X,L) \leqslant \mathfrak{R}_{\tau}^{+}(tP,X) + \frac{\text{Area}((t + \beta \tau ) P) - \text{Area}(t P)}{covol(L)} \textit{.} 
\end{equation}
Yet, one has the following two equations: 
\begin{equation}
\label{chap4:eq30}
\text{Area}((t - \beta \tau ) P) - \text{Area}(t P) = \text{Area}(P)(- 2 \beta \tau t + \beta^{2} \tau^{2}) 
\end{equation}
and
\begin{equation}
\label{chap4:eq31}
\text{Area}((t + \beta \tau ) P) - \text{Area}(t P) = \text{Area}(P)( 2 \beta \tau t + \beta^{2} \tau^{2}) \textit{.}
\end{equation}
Yet here $\tau = \frac{\log(t)^{\gamma}}{t}$ and so one gets that for every $T$ large enough, for every $ T \geqslant t \geqslant \alpha T$, 
\begin{equation}
\label{chap4:eq32}
| \frac{\pm 2 \beta \tau t + \beta^{2} \tau^{2}}{\sqrt{V(L^{\perp},t)}} | \leqslant M \log(T)^{\gamma - \frac{1}{2}}
\end{equation}
where $M> 0$. To obtain this last equation we have used the fact that $\sqrt{V(L^{\perp},t)} = \Theta(\sqrt{\log(t)})$ that is given by proposition $\ref{chap4:prop2}$. \\
By using Equation $\ref{chap4:eq29}$, Equation $\ref{chap4:eq30}$, Equation $\ref{chap4:eq31}$ and Equation $\ref{chap4:eq32}$, one gets the wanted result.
\end{proof}
The study is now reduced to the study, when $T \rightarrow \infty$, of the two quantities $\frac{\mathfrak{R}_{\tau}^{\pm}(tP,X)}{\sqrt{V(L^{\perp},t)}}$. In the next section, we are going to reduce the study of $\frac{\mathfrak{R}_{\tau}^{+}(tP,X)}{\sqrt{V(L^{\perp},t)}}$ to the study of the Fourier series $\frac{S(L^{\perp},X,t,T)}{\sqrt{V(L^{\perp},t)}}$ and the results are also going to hold for $\frac{\mathfrak{R}_{\tau}^{-}(tP,X)}{\sqrt{V(L^{\perp},t)}}$. 
%en fait, il faudrait faire le lemme conséquence (on se ramène à r+-) d'abord puis présenter le lemme encadrement après ... refaire la construction
\subsubsection{Simplification of the Fourier series $\frac{\mathfrak{R}_{\tau}^{+}(tP,X)}{\sqrt{V(L^{\perp},t)}}$}
%plus petit que t la norme ? ou T ? t log(t) ? T log(T) ? à voir ? et aussi le >0 ? faudra faire une remarque quand on passera à \Gamma\alpha parce qu'on prendra le k_{1} dans Z > 0, différent ! 
The main objective of this subsubsection is to prove the following proposition: 
\begin{prop}
\label{chap4:prop6}
Let us suppose that $L$ is admissible. Then, one has, for all $\gamma > 0$,
$$
 \mathbb{E}_{X \in \mathbb{R}^{2}/L}\left( (\frac{\mathfrak{R}_{\tau}^{+}(tP,X)-S(L^{\perp},t,X)}{\sqrt{V(t,L^{\perp})}})^{2} \right) \rightarrow 0
$$
where $\tau = \frac{\log(t)^{\gamma}}{t}$ and where the convergence towards $0$ is uniform in $t$ such that $\alpha T \leqslant t \leqslant T$ and when $T \rightarrow \infty$.
\end{prop}
To do this, we will need several intermediate lemmas. We will use tools from Fourier Analysis, integration calculus and Geometry of numbers. \\
$\emph{In the rest of this section, we are going to suppose that L is admissible.}$ \\
$\emph{Under this hypothesis, the dual lattice of L is also going to be admissible.}$

\paragraph{Reduction of $\mathfrak{R}_{\tau}^{+}$ to a finite sum.}
Let us introduce the following sum: 
\begin{equation}
\label{chap4:eq34}
S_{1}(L^{\perp},t,X) = \frac{1}{\pi^{2} \textit{covol}(L)} \sum_{\substack{l \in L^{\perp} \\ 0 < \lVert l \rVert \leqslant t  }} \frac{\hat{\omega}(\tau l) \sin(2 \pi l_{1} t^{+} ) \sin(2 \pi l_{2} t^{+})e^{2 i  \pi <l,X> }}{l_{1} l_{2}} 
\end{equation}
where $t^{+} = t + \beta \tau$. 
Then the main objective of this paragraph is to prove the following lemma: 
\begin{lemma}
\label{chap4:lemme4}
One has for all $T$ large enough,
$$
| \frac{|\mathfrak{R}_{\tau}^{+}(tP,X)-S_{1}(L^{\perp},t,X)|}{\sqrt{V(t,L^{\perp})}} | =  O(\log(T)^{\gamma - \frac{1}{2}}) \textit{.}
$$
where $O$ is uniform in $X \in \mathbb{R}^{2}$ and in $t$ when $\alpha T \leqslant t \leqslant T$.
\end{lemma}
Recall that $\gamma$ is a fixed parameter such that $ 0 < \gamma < \frac{1}{2}$. \\
To prove this lemma, we are going to need a lemma that gives an upper estimate of certain integrals: 
\begin{lemma}
\label{chap4:lemme5}
Let us set for $N \geqslant 2 $ and for $C > 0$,  
\begin{equation}
\label{chap4:eq35}
D_{N,C} = \{ x=(x_{1},x_{2}) \in (\mathbb{R}^{+})^{2} \text{ } | \text{ }  \lVert x \rVert \geqslant N \text{, } |\text{Num}(x)| \geqslant C  \} \textit{.}
\end{equation}
Then one has, for all $g: \mathbb{R}^{+} \longrightarrow \mathbb{R}^{+} $ measurable, there exists $M > 0$ such that for every $N$ large enough, 
$$\int_{x \in D_{N,C}} \frac{g(\lVert x \rVert)}{|x_{1}x_{2}|}dx_{1}dx_{2} \leqslant M \int_{r=N}^{\infty} g(r) \frac{\log(r)}{r} dr \textit{.} $$
\end{lemma}
\begin{proof}
Let us call: 
\begin{equation}
\label{chap4:eq36}
J(N,C) = \int_{x \in D_{N,C}} \frac{g(\lVert x \rVert)}{|x_{1}x_{2}|}dx_{1}dx_{2} \textit{.}
\end{equation}
By passing to polar coordinates, we get that for every $N$ large enough: 
\begin{equation}
\label{chap4:eq37}
J(N,C) \leqslant K_{1} \int_{r \geqslant N} \frac{g(r)}{r} (\int_{\arcsin(2 \frac{C}{r^{2}})}^{\frac{\pi}{2}} \frac{d \theta}{\sin(\theta)}) dr
\end{equation}
where $K_{1} > 0$.\\
Yet one calculates that:
\begin{equation}
\label{chap4:eq38}
(\int_{\arcsin(2  \frac{C}{r^{2}})}^{\frac{\pi}{2}} \frac{d \theta}{\sin(\theta)}) = \frac{1}{2}\log \big( \frac{(1+\sqrt{1-\frac{2  C}{r^{2}}})^{2}}{4  \frac{C}{r^{2}}} \big) \textit{.}
\end{equation}
With the Inequality ($\ref{chap4:eq37}$), we get that: 
\begin{equation}
\label{chap4:eq39}
J(N,C) \leqslant K_{2} \int_{r \geqslant N} \frac{g(r) \log(r) }{r} dr 
\end{equation}
where $K_{2} > 0$ and so we get the wanted result.
\end{proof}
We can now prove Lemma $\ref{chap4:lemme4}$: 
\begin{proof}[Proof of Lemma $\ref{chap4:lemme4}$]
$L$ is supposed to be admissible and so $L^{\perp}$ must also be admissible. \\
So there exists $C > 0$ such that for all $x \in L^{\perp}-\{ 0 \}$, 
\begin{equation}
\label{chap4:eq40}
| \text{Num}(x) | \geqslant C \textit{.}
\end{equation}
Furthermore, $\omega$ is supposed to be regular and of compact support. In consequence, $\hat{\omega}$ belongs to the Schwartz space and so for all $A > 2$, there exists $M > 0$ such that for all $x \in \mathbb{R}^{2}$: 
\begin{equation}
\label{chap4:eq41}
| \tilde{\omega}(x) | \leqslant \frac{M}{1 + \lVert x \rVert^{A}} \textit{.}
\end{equation}
Proposition $\ref{chap4:prop4}$ and the fact that $\mathbf{O}_{\tau}^{+} = (t+ \beta \tau) P$ with $P = \text{Rect}(a,b)$ give us that: 
\begin{equation}
\label{chap4:eq42}
|\mathfrak{R}_{\tau}^{+}(tP,X)-S_{1}(L^{\perp},t,X)| \leqslant K_{1} \sum_{\substack{l \in L^{\perp} \\ \lVert l \rVert > t}} \frac{1}{| l_{1} l_{2} | (1 + \lVert \tau l \rVert^{A})} 
\end{equation}
where $K_{1} > 0$. We have also used implicitly Equation $(\ref{chap4:eq23})$.  \\ %rajouter \lVert l \rVert > 0 et après E à enlever et changer lemme estimation (et donc prop au début ; juste mettre o de machin) 
Yet, one has also: 
\begin{equation}
\label{chap4:eq43}
 \sum_{\substack{l \in L^{\perp} \\ \lVert l \rVert > t}} \frac{1}{| l_{1} l_{2} | (1 + \lVert \tau l \rVert^{A})} \leqslant K_{2} \left( \int_{x \in D_{t,C}} \frac{1}{|x_{1}x_{2}|(1 + \lVert \tau x \rVert^{A})}dx_{1}dx_{2} + \int_{\substack{ \lVert x \rVert > t \\ |Num(x)| < C }} \frac{1}{C(1+\lVert \tau x \rVert^{A})} dx \right)
\end{equation}
where $K_{2} > 0$. \\
Lemma $\ref{chap4:lemme5}$ and a quick calculation give us then that for all $T$ large enough (we have to keep in mind that $\alpha T \leqslant t \leqslant T$): 
\begin{equation}
\label{chap4:eq44}
 \sum_{\substack{l \in L^{\perp} \\ \lVert l \rVert \leqslant t}} \frac{1}{| l_{1} l_{2} | (1 + \lVert \tau l \rVert^{A})} \leqslant K_{3} \left( \int_{r=t}^{\infty} \frac{\log(r)}{r (1 + (\tau r)^{A})} dr + \int_{r=t}^{\infty} \frac{1}{r^{2}(1 + (\tau r)^{A})} dr \right) 
\end{equation}
where $K_{3} > 0$. \\
However, one has also that:  
\begin{equation}
\label{chap4:eq45}
 \int_{r=t}^{\infty} \frac{\log(r)}{r (1 + (\tau r)^{A})} dr \leqslant \frac{K_{4} \log(t)}{ \tau^{A} t^{A}} 
\end{equation}
where $K_{4} > 0$ (and depends on $A$).\\
Yet one has that: $\tau = \frac{\log(t)^{\gamma}}{t}$ with $ \frac{1}{2}> \gamma > 0$. \\
By using Equation ($\ref{chap4:eq42}$), Equation ($\ref{chap4:eq44}$) and Equation ($\ref{chap4:eq45}$) and by using the fact that $V(L^{\perp},t) = \Theta(\log(t))$, one gets that: 
\begin{equation}
\label{chap4:eq46}
\frac{|\mathfrak{R}_{\tau}^{+}(tP,X)-S_{1}(L^{\perp},t,X)|}{\sqrt{V(t,L^{\perp})}} \leqslant K_{5} \left( \frac{\log(t)^{\frac{1}{2}}}{\log(t)^{ A \gamma}} + \frac{1}{\sqrt{\log(t)}} \right) 
\end{equation}
where $K_{5} > 0$ and depends on $A$. \\
Then, we can take $A$ large enough so one has that: 
\begin{equation}
\label{chap4:eq47}
\frac{\log(t)^{\frac{1}{2}}}{\log(t)^{ A \gamma}} \leqslant \frac{K_{6}}{\log(t)}  
\end{equation}
with $K_{6} > 0$.\\
By using Equation ($\ref{chap4:eq46}$) and Equation ($\ref{chap4:eq47}$), one has that: 
\begin{equation}
\label{chap4:eq48}
\frac{|\mathfrak{R}_{\tau}^{+}(tP,X)-S_{1}(L^{\perp},t,X)|}{\sqrt{V(t,L^{\perp})}}  \leqslant \frac{K_{7}}{\sqrt{\log(T)}}
\end{equation}
with $K_{7} > 0$ and for $T$ large enough. Thus the wanted result. 
\end{proof}
We are now brought back to the study of $S_{1}(L^{\perp},t,X)$ and now we are going, mainly, to "center" it on the prime vectors of $L^{\perp}$.
\paragraph{Centring of $S_{1}(L^{\perp},t,X)$ over prime vectors.} %reprendre ici
Before stating the main lemma of this paragraph, let us make a small remark. Because $L^{\perp}$ is admissible and because $\tilde{\omega}$ is spherically symmetric, by parity, one has that: 
\begin{equation}
\label{chap4:eq49}
S_{1}(L^{\perp},t,X) = \frac{2}{\pi^{2} \textit{covol}(L)} \sum_{l \in J_{1}(L^{\perp},t) } \frac{\tilde{\omega}(\tau l) \sin(2 \pi l_{1} t^{+} ) \sin(2 \pi l_{2} t^{+}) \cos(2   \pi <l,X>) }{l_{1} l_{2}} 
\end{equation}
where 
\begin{equation}
\label{chap4:eq50}
J_{1}(L^{\perp},t) = \{ l \in L^{\perp} \text{ | } 0 < \lVert l \rVert \leqslant t \text{ and } l_{1} > 0 \} \textit{.}
\end{equation}
Let us define: 
\begin{equation}
\label{chap4:eq51}
S_{2}(L^{\perp},t,X) = \frac{2}{\pi^{2} \textit{covol}(L)} \sum_{l \in J_{2}(L^{\perp},t) } \frac{Z_{l}(X,t)}{l_{1} l_{2}} 
\end{equation}
where
\begin{equation}
\label{chap4:eq52}
Z_{l}(X,t) = \sum_{k=1}^{\infty} \frac{\tilde{\omega}(k \tau l) \sin(2 k \pi l_{1} t^{+} ) \sin(2 k \pi l_{2} t^{+}) \cos(2  k \pi <l,X>)}{k^{2}} 
\end{equation}
and we recall that
\begin{equation*}
J_{2}(L^{\perp},t) = \{ l \in L^{\perp} \text{ | } 0 < \lVert l \rVert \leqslant t \text{ , } l \text{ prime} \text{ and } l_{1} > 0 \} \textit{.}
\end{equation*}
Then the main statement of this subsection is the following lemma: 
\begin{lemma}
\label{chap4:lemme6}
For every $T$ large enough, one has: 
$$
 |\frac{|S_{2}(L^{\perp},t,X)-S_{1}(L^{\perp},t,X)|}{\sqrt{V(L^{\perp},t)}}| = O(\log(T)^{\gamma - \frac{1}{2}}) 
$$
with $O$ being uniform in $X \in \mathbb{R}^{2}$ and in $t$ such that $ \alpha T \leqslant t \leqslant T$. 
\end{lemma}
It basically says that the essential information of $S_{1}$ is contained in its prime terms.
\begin{proof}
One has: 
\begin{equation}
\label{chap4:eq54}
|S_{2}(L^{\perp},t,X)-S_{1}(L^{\perp},t,X)| \leqslant K_{1} \sum_{ \substack{ l \in L^{\perp} \\ \lVert l \rVert \geqslant t }} \frac{| \tilde{\omega}(\tau l) |}{| l_{1} l_{2} |} 
\end{equation}
where $K_{1} > 0$. \\
So, we can apply the same argument as before in the proof of Lemma $\ref{chap4:lemme4}$ and one gets that: 
\begin{equation}
\label{chap4:eq56}
 |\frac{|S_{2}(L^{\perp},t,X)-S_{1}(L^{\perp},t,X)|}{\sqrt{V(L^{\perp},t)}}| = O(\log(T)^{\gamma - \frac{1}{2}}) \textit{.}
\end{equation}
\end{proof}
So, the study of $S_{1}$, when $T \rightarrow \infty$, is now reduced to the study of $S_{2}$, when $T \rightarrow \infty$. In the next paragraph, we are going to simplify the sum $S_{2}$.
\paragraph{Replacing $\tilde{\omega}$ by $1$ in $S_{2}$.}
Let us define: 
\begin{equation}
\label{chap4:eq57}
S_{3}(L^{\perp},t,X) = \frac{2}{\pi^{2} \textit{covol}(L)} \sum_{l \in J_{2}(L^{\perp},t) } \frac{\tilde{Z}_{l}(X,t)}{l_{1} l_{2}} 
\end{equation}
where
\begin{equation}
\label{chap4:eq58}
\tilde{Z}_{l}(X,t) = \sum_{k=1}^{\infty} \frac{\sin(2 k \pi l_{1} t^{+} ) \sin(2 k \pi l_{2} t^{+}) \cos(2  k \pi <l,X>)}{k^{2}}
\end{equation} 
(we have replaced $\tilde{\omega}$ by its value at $0$, which is $1$ because $\int_{\mathbb{R}^{2}} \omega = 1$). \\
The main statement of this subsection is the following lemma: 
\begin{lemma}
\label{chap4:lemme7}
Let us set: 
\begin{equation}
\label{chap4:eq59}
\Delta(L^{\perp},t) = \int_{X \in \mathbb{R}^{2}/L^{\perp}} (\frac{|S_{2}(L^{\perp},t,X)-S_{3}(L^{\perp},t,X)|}{\sqrt{V(L^{\perp},t)}})^{2} d\tilde{\lambda}_{2}(X) \textit{.}
\end{equation}
Then $\Delta(L^{\perp},t)$ goes to $0$ uniformly in $t$ that are such that $ \alpha T \leqslant t \leqslant T$ with $T \rightarrow \infty$. 
\end{lemma}
\begin{proof}
The Parseval formula applies here and gives us that: 
\begin{equation}
\label{chap4:eq60}
\Delta(L^{\perp},t) \leqslant \frac{K_{1}}{\log(t)} \sum_{l \in J_{1}(L^{\perp},t)} (\frac{\tilde{\omega}(\tau l) - 1}{l_{1}l_{2}})^{2}  
\end{equation}
where $K_{1} > 0$ and we have used the fact that $V(L^{\perp},t) = \Theta(\log(t))$. \\
Let us take $\tilde{\gamma}$ such that $ \gamma < \tilde{\gamma} < \frac{1}{2}$ and let us set $t_{1} = \frac{t}{\log(t)^{\tilde{\gamma}}}$. \\
Then, with this notation, one has, thanks to Equation $(\ref{chap4:eq60})$: 
\begin{equation}
\label{chap4:eq61}
\Delta(L^{\perp},t) \leqslant K_{2}( \Delta_{1} + \Delta_{2}) 
\end{equation}
where $K_{2} > 0$, 
\begin{equation}
\label{chap4:eq62}
\Delta_{1} = \frac{1}{\log(t)} \sum_{l \in J_{1}(L^{\perp},t_{1})} (\frac{\tau \lVert l \rVert }{l_{1}l_{2}})^{2}  \textit{ and}
\end{equation}

\begin{equation}
\label{chap4:eq63}
\Delta_{2} = \frac{1}{\log(t)} \sum_{l \in J_{1}(L^{\perp},t) - J_{1}(L^{\perp},t_{1})} \frac{1}{l_{1}^{2} l_{2}^{2}} \textit{.}
\end{equation}
Then Lemma $\ref{chap4:lemme2}$, the fact $\tau = \frac{\log(t)^{\gamma}}{t}$ and the definition of $t_{1}$ give us that for $T$ large enough (and for $t$ such that $ \alpha T \leqslant t \leqslant T $): 
\begin{equation}
\label{chap4:eq64}
\Delta_{1} \leqslant K_{3} \log(t)^{2(\gamma - \tilde{\gamma})} 
\end{equation}
with $K_{3}> 0$ and 
\begin{equation}
\label{chap4:eq65}
\Delta_{2} \leqslant \frac{K_{4}}{\log(t)} \int_{ t \log(t)^{- \tilde{\gamma}}}^{t} \frac{1}{r}dr = O( \frac{1}{\sqrt{\log(t)}}) \textit{.}
\end{equation}
with $K_{4} > 0$. \\
So, one gets finally that, when $T \rightarrow \infty$, 
$$\Delta_{1} \rightarrow 0 \textit{ and } \Delta_{2} \rightarrow 0 $$
uniformly in $t$ for $\alpha T \leqslant t \leqslant T$. \\
\end{proof}
So, now, we are brought back to the study of $S_{3}$ when $T \rightarrow \infty$. \\
In the next paragraph, we are going to simply $S_{3}$ and its study will be reduced to the study of $S$ when $T \rightarrow \infty$ as wanted.
\paragraph{Replacing $t^{+}$ by $t$ and proof of proposition $\ref{chap4:prop6}$.} 
We recall that we have defined $S$ by the Equation ($\ref{chap4:eq19}$). The main objective of this paragraph is to prove the following lemma: 
\begin{lemma}
\label{chap4:lemme8}
Let us call (here) $\Delta$ the following quantity: 
\begin{equation}
\label{chap4:eq66}
\Delta(L^{\perp},t) = \sqrt{\int_{X \in \mathbb{R}^{2}/L^{\perp}} (\frac{|S_{3}(L^{\perp},t,X)-S(L^{\perp},t,X)|}{\sqrt{V(L^{\perp},t)}})^{2} d\tilde{\lambda}_{2}(X)} \textit{.}
\end{equation}
Then $\Delta(L^{\perp},t)$ goes to $0$ uniformly in $t$ that are such that $ \alpha T \leqslant t \leqslant T$ with $T \rightarrow \infty$. 
\end{lemma}
Mainly, this lemma says that we can replace $t^{+}$ in the two $\sin$ terms by $t$. It is reasonable insofar as $t^{+}$ is very close to $t$. 
\begin{proof}
Let us set: 
\begin{equation}
\label{chap4:eq67}
S_{4}(L^{\perp},t,X) = \frac{2}{\pi^{2} \textit{covol}(L)} \sum_{l \in J_{2}(L^{\perp},t) } \frac{W_{l}(X,t)}{l_{1} l_{2}} 
\end{equation}
where 
\begin{equation}
\label{chap4:eq68}
W_{l}(X,t) = \sum_{k=1}^{\infty} \frac{\sin(2 k \pi l_{1} t ) \sin(2 k \pi l_{2} t^{+}) \cos(2  k \pi <l,X>)}{k^{2}} \textit{.}
\end{equation} 
Then the triangle inequality gives us that: 
\begin{equation}
\label{chap4:eq69}
\Delta \leqslant \sqrt{\Delta_{1}} + \sqrt{\Delta_{2}}
\end{equation}
where 
\begin{equation}
\label{chap4:eq70}
\Delta_{1} = \int_{X \in \mathbb{R}^{2}/L^{\perp}} (\frac{|S_{3}(L^{\perp},t,X)-S_{4}(L^{\perp},t,X)|}{\sqrt{V(L^{\perp},t)}})^{2} d\tilde{\lambda}_{2}(X)
\end{equation}
and
\begin{equation}
\label{chap4:eq71}
\Delta_{2} = \int_{X \in \mathbb{R}^{2}/L^{\perp}} (\frac{|S_{4}(L^{\perp},t,X)-S(L^{\perp},t,X)|}{\sqrt{V(L^{\perp},t)}})^{2} d\tilde{\lambda}_{2}(X) \textit{.}
\end{equation}
Let us take $\tilde{\gamma}$ such that $ \gamma < \tilde{\gamma} < \frac{1}{2}$ and let us set $t_{1} = \frac{t}{\log(t)^{\tilde{\gamma}}}$. \\
The Parseval formula and the mean value theorem apply here and give us that: 
\begin{equation}
\label{chap4:eq72}
\Delta_{1} \leqslant K_{2}( \Delta_{1,1} + \Delta_{1,2}) 
\end{equation}
where $K_{2} > 0$, 
\begin{equation}
\label{chap4:eq73}
\Delta_{1,1} = \frac{1}{\log(t)} \sum_{l \in J_{1}(L^{\perp},t_{1})} (\frac{\tau \lVert l \rVert }{l_{1}l_{2}})^{2}  \textit{, }
\end{equation}

\begin{equation}
\label{chap4:eq74}
\Delta_{1,2} = \frac{1}{\log(t)} \sum_{l \in J_{1}(L^{\perp},t) - J_{1}(L^{\perp},t_{1})} \frac{1}{l_{1}^{2} l_{2}^{2}} 
\end{equation}
and
\begin{equation}
\label{chap4:eq75}
\Delta_{2} \leqslant K_{2}( \Delta_{1,1} + \Delta_{1,2}) \textit{.} 
\end{equation}
We have used the fact that $V(L^{\perp},) = \Theta(\log(t))$. \\
Then Equation $(\ref{chap4:eq64})$, Equation $(\ref{chap4:eq65})$, Equation ($\ref{chap4:eq69}$), Equation $(\ref{chap4:eq75})$ and Equation ($\ref{chap4:eq72}$) give that $\Delta$ goes uniformly to $0$ in $t$ for $t$ such that $\alpha T \leqslant t \leqslant T$
\end{proof}
We can now prove Proposition $\ref{chap4:prop6}$.
\begin{proof}[Proof of proposition $\ref{chap4:prop6}$]
Proposition $\ref{chap4:prop6}$ is a direct consequence of the lemmas $\ref{chap4:lemme4}$, $\ref{chap4:lemme6}$, $\ref{chap4:lemme7}$ and $\ref{chap4:lemme8}$. 
\end{proof}
\subsubsection{Case of $t^{-}$ and proof of Proposition $\ref{chap4:prop3}$}
By following exactly the same approach that has been used to prove proposition $\ref{chap4:prop6}$, we prove the following proposition: 
\begin{prop}
\label{chap4:prop7}
Let us suppose that $L$ is admissible. Then, one has, for all $\gamma > 0$,
$$
\mathbb{E}_{X \in \mathbb{R}^{2}/L} \left( ( \frac{|\mathfrak{R}_{\tau}^{-}(tP,X)-S(L^{\perp},t,X)|}{\sqrt{V(t,L^{\perp})}})^{2} \right) \rightarrow 0 \textit{.}
$$
where $\tau = \frac{\log(t)^{\gamma}}{t}$ and where the convergence towards $0$ is uniform in $t$ such that $\alpha T \leqslant t \leqslant T$ and when $T \rightarrow \infty$.
\end{prop}
We are now able to prove Proposition $\ref{chap4:prop3}$: 
\begin{proof}[Proof of Proposition $\ref{chap4:prop3}$]
It is a direct consequence of Proposition $\ref{chap4:prop5}$, Proposition $\ref{chap4:prop6}$ and Proposition $\ref{chap4:prop7}$ and of the triangle inequality.
\end{proof} 
So, finally, it is enough to study the behaviour, when $T \rightarrow \infty$, of $S$. It is the object of the next subsection.  %must be modified, theorem 1 laisser celui-là ?
%ou un autre ? indépendance ?
\subsection{Asymptotic behaviour of $S$}
According to the proof of proposition $\ref{chap4:prop3}$, by replacing $L^{\perp}$ by $L$ (if one is admissible, the other also is and conversely), to show the second assertion of Theorem $\ref{chap4:thm1}$, one only needs to prove that 
$\mathbb{E}_{X \in \mathbb{R}^{2}/L^{\perp}}\left( (\frac{S(L,X,t)}{\sqrt{V(L,t)}})^{2} \right)$ converges in distribution and in probability towards  $\frac{1}{4 \pi^{4} \text{Covol}(L^{\perp})^{2}}$. Furthermore, this limit constant will be the same for $\mathbb{E}_{X \in \mathbb{R}^{2}/L^{\perp}}\left( (\frac{S(L,X,t)}{\sqrt{V(L,t)}})^{2} \right)$ and $\mathbb{E}_{X \in \mathbb{R}^{2}/L^{\perp}}\left( (\frac{S(L,X,t)}{\sqrt{V(L,t)}})^{2} \right)$. \\%il faudra reformuler toutes les propositions avant et changer légèrement la fin des preuves
With the Parseval formula, we see that we only need to prove the following proposition: 
\begin{prop}
\label{chap4:prop8}
$$ G(L,t) = \frac{2}{\pi^{4} \text{Covol}(L^{\perp})^{2} } \sum_{l \in J_{2}(L,t)} \frac{1}{(l_{1}l_{2})^{2}} \sum_{k=1}^{\infty} \frac{(\sin(2 \pi k t l_{1}) \sin(2 \pi k t l_{2}))^{2} }{k^{4}} $$ normalized by $V(L,t)$ converges in distribution and in probability towards $\frac{1}{4 \pi^{4} \text{Covol}(L^{\perp})^{2}}$.
\end{prop}
The rest of this subsection is dedicated to proving this last proposition. 
\subsubsection{A small remark and approach}
By using the fact that one has for every $x \in \mathbb{R}$,  
$$\sin(x)^{2} = \frac{1 - \cos(2 x) }{2} \textit{, } $$
one gets that: 
\begin{equation}
\label{chap4:eq90}
G(L,t) = G_{1}(L,t) - G_{2}(L,t) - G_{3}(L,t) + G_{4}(L,t) 
\end{equation}
where 
\begin{equation}
\label{chap4:eq91}
G_{1}(L,t) = \frac{1}{2 \pi^{4}\text{Covol}(L^{\perp})^{2}} \sum_{l \in J_{2}(L,t)} \frac{1}{(l_{1}l_{2})^{2}} \sum_{k=1}^{\infty} \frac{1 }{k^{4}} \textit{, }
\end{equation} 
\begin{equation}
\label{chap4:eq92}
G_{2}(L,t) = \frac{1}{2 \pi^{4}\text{Covol}(L^{\perp})^{2}} \sum_{l \in J_{2}(L,t)} \frac{1}{(l_{1}l_{2})^{2}} \sum_{k=1}^{\infty} \frac{\cos (4 \pi k t l_{1}) }{k^{4}} \textit{, }
\end{equation} 
\begin{equation}
\label{chap4:eq93}
G_{3}(L,t) = \frac{1}{2 \pi^{4}\text{Covol}(L^{\perp})^{2}} \sum_{l \in J_{2}(L,t)} \frac{1}{(l_{1}l_{2})^{2}} \sum_{k=1}^{\infty} \frac{\cos (4 \pi k t l_{2}) }{k^{4}} \textit{ and }
\end{equation} 
\begin{equation}
\label{chap4:eq94}
G_{4}(L,t) = \frac{1}{2 \pi^{4}\text{Covol}(L^{\perp})^{2}} \sum_{l \in J_{2}(L,t)} \frac{1}{(l_{1}l_{2})^{2}} \sum_{k=1}^{\infty} \frac{\cos (4 \pi k t l_{1})\cos (4 \pi k t l_{2}) }{k^{4}} \textit{. }
\end{equation} 
To get the validity of Proposition $\ref{chap4:prop8}$, we only need to show the three following propositions: 
\begin{prop}
\label{chap4:prop9}
$$\frac{G_{1}(L,t)}{V(L,t)} \rightarrow \frac{1}{4 \pi^{4}\text{Covol}(L^{\perp})^{2}} $$ when $t \rightarrow \infty$. 
\end{prop}
\begin{prop}
\label{chap4:prop10}
$ \frac{G_{2}(L,t)}{V(L,t)} $ and $ \frac{G_{3}(L,t)}{V(L,t)} $ converge in distribution and in probability towards $0$ when $t$ is distributed according to the probability measure $\frac{1}{T} \rho(\frac{t}{T}) dt$ and when $T \rightarrow \infty$.
\end{prop}
\begin{prop}
\label{chap4:prop11}
$ \frac{G_{4}(L,t)}{V(t,T)} $ converge in distribution and in probability towards $0$ when $t$ is distributed according to the probability measure $\frac{1}{T} \rho(\frac{t}{T}) dt$ and when $T \rightarrow \infty$.
\end{prop}
Proposition $\ref{chap4:prop9}$ is just a use of different definitions whereas Proposition $\ref{chap4:prop10}$ and Proposition $\ref{chap4:prop11}$ work because, basically, there are, relatively to $G_{1}$, additional oscillatory terms that make the normalized sum go to $0$. \\
Now we are going to prove in this order these three propositions.
\subsubsection{Proof of Proposition $\ref{chap4:prop9}$}
The proof of this proposition is straightforward: 
\begin{proof}[Proof of Proposition $\ref{chap4:prop9}$]
One has: 
$$
G_{1}(L,t) = \frac{1}{2 \pi^{4} \text{Covol}(L^{\perp})^{2}} \sum_{l \in J_{2}(L,t)} \frac{1}{(l_{1}l_{2})^{2}} \sum_{k=1}^{\infty} \frac{1 }{k^{4}}
$$
where 
$$ J_{2}(L,t) = \{ l \in L \text{ } | \text{ } 0 < \lVert l \rVert \leqslant t \text{ , } l \text{ prime} \text{ and } l_{1} > 0 \} $$
and 
$$ V(L,t)= \sum_{\substack{ l \in L \\ 0 < \lVert l \rVert \leqslant t }} \frac{1}{l_{1}^{2} l_{2}^{2}} \textit{.} $$
By making $t$ goes to infinity, one gets the wanted result.
\end{proof} 
Next, we are going to prove Proposition $\ref{chap4:prop10}$.
\subsubsection{Proof of Proposition $\ref{chap4:prop10}$}
%on a besoin du lemme calculatoire (comme un peu lemme 2) normal puis 
%faire la preuve du premier lemme préparatoire 
%puis direct preuve prop 10, faire un seul des deux cas
%rédiger normal et peut-être rajouter un lemme calculatoire avant 
Before giving the proof of Proposition $\ref{chap4:prop10}$, we need several lemmas. \\
Let us introduce:
$$ \tilde{G}_{2}(L,t,T) = \frac{1}{2 \pi^{4}\text{Covol}(L^{\perp})^{2}} \sum_{l \in J_{2}(L,T)} \frac{1}{(l_{1}l_{2})^{2}} \sum_{k=1}^{\infty} \frac{\cos (4 \pi k t l_{1}) }{k^{4}} $$
and 
$$ \tilde{G}_{3}(L,t,T) = \frac{1}{2 \pi^{4}\text{Covol}(L^{\perp})^{2}} \sum_{l \in J_{2}(L,T)} \frac{1}{(l_{1}l_{2})^{2}} \sum_{k=1}^{\infty} \frac{\cos (4 \pi k t l_{2}) }{k^{4}} \textit{.} $$
 The first lemma basically says that we can study $\frac{\tilde{G}_{2}(L,t,T)}{V(L,t)}$ and $\frac{\tilde{G}_{3}(L,t,T)}{V(L,t)}$ instead of studying, respectively, $\frac{G_{2}(L,t)}{V(L,t)}$ and $\frac{G_{3}(L,t)}{V(L,t)}$: 
\begin{lemma}
\label{chap4:lemme10}
One has, when $T \rightarrow \infty$, 
$$\mathbb{E} \left( |\frac{\tilde{G}_{2}(L,t,T)}{V(L,t)}-\frac{G_{2}(L,t)}{V(L,t)}| \right) \rightarrow 0$$ 
and 
$$\mathbb{E} \left( |\frac{\tilde{G}_{3}(L,t,T)}{V(L,t)}-\frac{G_{3}(L,t)}{V(L,t)}| \right) \rightarrow 0  $$
(Note that the expectation are calculated relatively to $t$ with $t$ being distributed according to the probability measure $\frac{1}{T} \rho(\frac{t}{T}) dt$).
\end{lemma}
\begin{proof}
We are only going to prove the first fact because the proof is going to be symmetrical in $(l_{1},l_{2})$. \\
One has that:
\begin{equation}
\label{chap4:eq101}
\frac{\tilde{G}_{2}(L,t,T)}{V(L,t)}-\frac{G_{2}(L,t)}{V(L,t)} = \frac{D}{V(L,t)} \sum_{l \in J_{2}(L,T) -J_{2}(L,t) } \frac{1}{(l_{1}l_{2})^{2}} \sum_{k=1}^{\infty} \frac{\cos (4 \pi k t l_{1}) }{k^{4}}
\end{equation}
where $D$ is a positive constant. \\
By integrating and because $t$ belongs to $[\alpha T, T]$, one gets that: 
\begin{equation}
\label{chap4:eq102}
\mathbb{E}\left( |\frac{\tilde{G}_{2}(L,t,T)}{V(L,t)}-\frac{G_{2}(L,t)}{V(L,t)}| \right) \leqslant \frac{D}{\log(T)} \mathbb{E} \left( \sum_{l \in J_{2}(L,T) -J_{2}(L,t) } \frac{1}{(l_{1}l_{2})^{2}} \sum_{k=1}^{\infty} \frac{1}{k^{4}} \right)
\end{equation}
where $D$ is a positive constant, possibly different from the previous one. We have also used Proposition $\ref{chap4:prop2}$. \\
Yet, one has that (see Lemma $\ref{chap4:lemme2}$):  
\begin{equation}
\label{chap4:eq103}
\sum_{l \in J_{2}(L,T) -J_{2}(L,t) } \frac{1}{(l_{1}l_{2})^{2}} \sum_{k=1}^{\infty} \frac{1}{k^{4}} \leqslant D \log(\frac{T}{t}) \textit{.}
\end{equation}
So, with Equation ($\ref{chap4:eq102}$), one has that: 
\begin{equation}
\label{chap4:eq104}
\mathbb{E}\left( |\frac{\tilde{G}_{2}(L,t,T)}{V(L,t)}-\frac{G_{2}(L,t)}{V(L,t)}| \right) \leqslant \frac{D}{\log(T)} \mathbb{E}( \log(\frac{T}{t})) \textit{.}
\end{equation}
Yet, a quick calculation gives us that: 
\begin{equation}
\label{chap4:eq105}
\mathbb{E}( \log(\frac{T}{t})) = O(1) 
\end{equation}
when $T \rightarrow \infty$. \\
So, thanks to Equation $(\ref{chap4:eq104})$ and thanks to Equation $(\ref{chap4:eq105})$, one gets the first wanted result. 
\end{proof}
The second lemma is an estimating one. 
\begin{lemma}
\label{chap4:lemme9}
For every $A > 0$, for every $C > 0$, one has: 
$$ \int_{ \substack{ A \leqslant \lVert l \rVert \leqslant T \\ l_{1} > 0,\textit{ } l_{2} > 0 \\ l_{1}l_{2} \geqslant C}} \frac{1}{l_{1}^{3} l_{2}^{2}} dl_{1} dl_{2} = O(T) \textit{ and } $$
$$ \int_{ \substack{ A \leqslant \lVert l \rVert \leqslant T \\ l_{1} > 0,\textit{ } l_{2} > 0 \\ l_{1}l_{2} \geqslant C }} \frac{1}{l_{1}^{2} l_{2}^{3}} dl_{1} dl_{2} = O(T) \textit{.} $$
\end{lemma}
\textit{Remark.} Thanks to this lemma, we can show very quickly that the expectation of $\frac{\tilde{G}_{2}(L,t,T)}{V(L,t)}$ and of $\frac{\tilde{G}_{2}(L,t,T)}{V(L,t)}$ tends to $0$ when $T \rightarrow \infty$.
\begin{proof}
By symmetry, one has: 
\begin{equation}
\label{chap4:eq95}
\int_{ \substack{ A \leqslant \lVert l \rVert \leqslant T \\ l_{1} > 0,\textit{ } l_{2} > 0 \\ l_{1}l_{2} \geqslant C}} \frac{1}{l_{1}^{3} l_{2}^{2}} dl_{1} dl_{2} = \int_{ \substack{ A \leqslant \lVert l \rVert \leqslant T \\ l_{1} > 0,\textit{ } l_{2} > 0 \\ l_{1}l_{2} \geqslant C }} \frac{1}{l_{1}^{2} l_{2}^{3}} dl_{1} dl_{2}
\end{equation}
So, we only have to prove the first equality of lemma $\ref{chap4:lemme9}$. \\
By passing into polar coordinates $(r, \theta)$, one has that: 
\begin{equation}
\label{chap4:eq96}
\int_{ \substack{ A \leqslant \lVert l \rVert \leqslant T \\ l_{1} > 0,\textit{ } l_{2} > 0 \\ l_{1}l_{2} \geqslant C}} \frac{1}{l_{1}^{3} l_{2}^{2}} dl_{1} dl_{2} = \int_{ \substack{ A \leqslant r \leqslant T \\ \frac{\pi}{2} > \theta > 0 \\ \sin(2 \theta) \geqslant \frac{2 C}{r^{2}} }} \frac{1}{r^{4}} \frac{1}{\cos^{3}(\theta) \sin^{2}(\theta)} dr d\theta \textit{.} 
\end{equation}
We see that, when $t \rightarrow \infty$, there are $\textit{a priori}$ two essential parts of this last integral: the first one is when $r$ is large and $\theta$ is closed to $\frac{1}{2} \arcsin( \frac{2 C}{r^{2}})$ ; the second is when $r$ is large and $\theta$ is closed to $\frac{\pi}{2} - \frac{1}{2} \arcsin( \frac{2 C}{r^{2}})$. \\
By using the facts that when $\theta \rightarrow 0$, $\sin(\theta) \sim \theta$ and when $\theta \rightarrow \frac{\pi}{2}$, $\cos(\theta) \sim \frac{\pi}{2}-\theta$ and by calculating, one gets that the first essential part, let us call it $A_{1}(T)$, is estimated as followed
\begin{equation}
\label{chap4:eq97}
A_{1}(T) = O(1)
\end{equation}
whereas the second essential part, let us call it $A_{2}(T)$, is estimated as followed  
\begin{equation}
\label{chap4:eq98}
A_{2}(T) = O(T) \textit{.}
\end{equation}
By using Equation ($\ref{chap4:eq96}$), one gets that: 
\begin{equation}
\label{chap4:eq99}
\int_{ \substack{ A \leqslant \lVert l \rVert \leqslant T \\ l_{1} > 0,\textit{ } l_{2} > 0 \\ l_{1}l_{2} \geqslant C}} \frac{1}{l_{1}^{3} l_{2}^{2}} dl_{1} dl_{2}  = O(T) \textit{.}
\end{equation}
%Continuez: coordonnées polaires et t rightarrow infty

\end{proof}
%Non en fait ce lemme ne sert à rien, il est trop brutale, il faut reprendre. ce n'est pas log(t) mais t
We can now tackle the proof of Proposition $\ref{chap4:prop10}$.
\begin{proof}[Proof of Proposition $\ref{chap4:prop10}$]
The proof will be symmetrical relatively to the transformation $l_{1} \leftarrow l_{2}$. So, we only need to give the proof of the result that concerns $G_{2}(L,t)$. \\
According to Lemma $\ref{chap4:lemme10}$ and Markov's inequality, we only need to see that $\frac{\tilde{G}_{2}(L,t,T)}{V(L,t)}$ converges in distribution and in probability towards $0$. \\
To obtain the fact that $\frac{\tilde{G}_{2}(L,t,T)}{V(L,t)} \rightarrow 0$ in probability, we are going to show that its second moment goes to $0$. \\
One has, according to the definition of $\frac{\tilde{G}_{2}(L,t,T)}{V(L,t)}$, that: 
\begin{equation}
\label{chap4:eq107}
\mathbb{E}  \big( (\frac{\tilde{G}_{2}(L,t,T)}{V(L,t)})^{2} \big) \leqslant D \mathbb{E}\left( \frac{1}{V(L,t)^{2}} \sum_{l,l' \in J_{2}(L,T)} \frac{1}{(l_{1}l_{2})^{2}}\frac{1}{(l'_{1}l'_{2})^{2}} \sum_{k,k' \geqslant 1} \frac{\cos (4 \pi k t l_{1}) }{k^{4}} \frac{\cos (4 \pi k' t l'_{1}) }{k'^{4}} \right) 
\end{equation}
where $D > 0$.\\
By integrating, by using a usual trigonometric formula and by using Proposition $\ref{chap4:prop2}$, one gets that: 
\begin{equation}
\label{chap4:eq108}
\mathbb{E}  \big( (\frac{\tilde{G}_{2}(L,t,T)}{V(L,t))})^{2} \big) \leqslant  O(A_{1}(T) + A_{2}(T)) \textit{.}
\end{equation}
where 
\begin{equation}
\label{chap4:eq109}
A_{1}(T) =  \sum_{l,l' \in J_{2}(L,T)} \frac{1}{(l_{1}l_{2})^{2}}\frac{1}{(l'_{1}l'_{2})^{2}} \sum_{k,k' \geqslant 1} \frac{1}{(k k')^{4}} \max(\frac{1}{\log(T)^{2}}, \frac{1}{\log(T)^{2}T|kl_{1}+k'l'_{1}| })
\end{equation}
and
\begin{equation}
\label{chap4:eq110}
A_{2}(T) =  \sum_{l,l' \in J_{2}(L,T)} \frac{1}{(l_{1}l_{2})^{2}}\frac{1}{(l'_{1}l'_{2})^{2}} \sum_{k,k' \geqslant 1} \frac{1}{(k k')^{4}} \max(\frac{1}{\log(T)^{2}}, \frac{1}{\log(T)^{2}T|kl_{1}-k'l'_{1}| }) \textit{.}
\end{equation}
To get the wanted result, it is enough to show that $A_{1}(T)$ and $A_{2}(T)$ tend to $0$ when $T \rightarrow \infty$. Furthermore, we are only going to show that $A_{2}(T)$ tend to $0$ when $T \rightarrow \infty$, the proof for $A_{1}(T)$ being symmetrical. \\
Let us remark that one has: 
\begin{equation}
\label{chap4:eq111}
A_{2}(T) \leqslant A_{2,1}(T) + A_{2,2}(T)
\end{equation}
where
\begin{equation}
\label{chap4:eq112}
A_{2,1}(T) = \sum_{l,l' \in J_{2}(L,T)} \frac{1}{(l_{1}l_{2})^{2}}\frac{1}{(l'_{1}l'_{2})^{2}} \sum_{\substack{ k,k' \geqslant 1 \\ \min(k,k') \geqslant \lceil \log(T) \rceil }} \frac{1}{(k k')^{4}} \max(\frac{1}{\log(T)^{2}}, \frac{1}{\log(T)^{2}T|kl_{1}-k'l'_{1}| })
\end{equation}
and
\begin{equation}
\label{chap4:eq113}
A_{2,2}(T) = \sum_{l,l' \in J_{2}(L,T)} \frac{1}{(l_{1}l_{2})^{2}}\frac{1}{(l'_{1}l'_{2})^{2}} \sum_{\lceil \log(T) \rceil \geqslant k,k' \geqslant 1 } \frac{1}{(k k')^{4}} \max(\frac{1}{\log(T)^{2}}, \frac{1}{\log(T)^{2}T|kl_{1}-k'l'_{1}| }) \textit{.}
\end{equation}
From Lemma $\ref{chap4:lemme2}$ and by using a usual equivalent, one gets that: 
\begin{equation}
\label{chap4:eq114}
A_{2,1}(T) \rightarrow 0
\end{equation}
when $T \rightarrow \infty$. As a consequence, we only need to look at the behaviour of $A_{2,2}(T)$. \\
There are two types of terms in $A_{2,2}(T)$: those such that $k l_{1}$ is close to $k' l'_{1}$, for example at a distance less than $\frac{1}{T}$ and the others. \\
In the first case, it forms a sum that is estimated by  $$\frac{D}{\log(T)^{2}} \int_{\substack{l_{1} > 0 \textit{, } l_{2} > 0 \\ l_{1} l_{2} \geqslant C \\  D \leqslant \lVert l \rVert \leqslant T \log(T)}} \frac{1}{l_{1}^{2}l_{2}^{2}} \frac{1}{T l_{1}}dl_{1}dl_{2}$$ with $D > 0$ and this last quantity goes to zero according to Lemma $\ref{chap4:lemme9}$ (we have first integrated over $l'_{2}$ and then used the fact $k l_{1}$ is close to $k' l'_{1}$).  \\
So, finally, we only need to show that the following quantity goes to $0$ when $T$ goes to infinity: 
\begin{equation}
\label{chap4:eq115}
J(T) = \sum_{\substack{ \lceil \log(T) \rceil \geqslant k,k' \geqslant 1 \\ |kl_{1} - k'l'_{1}|\geqslant \frac{1}{T} }} \frac{1}{(l_{1}l_{2})^{2}}\frac{1}{(l'_{1}l'_{2})^{2}} \frac{1}{(k k')^{4}} \max(\frac{1}{\log(T)^{2}}, \frac{1}{\log(T)^{2}T|kl_{1}-k'l'_{1}| }) \textit{.}
\end{equation}
Yet, one has that: 
\begin{equation}
\label{chap4:eq116}
J(T) \leqslant J_{1}(T) + J_{2}(T)
\end{equation}
where 
\begin{equation}
\label{chap4:eq117}
J_{1}(T)= \sum_{\substack{ \lceil \log(T) \rceil \geqslant k,k' \geqslant 1 \\ |kl_{1} - k'l'_{1}|\geqslant 1 }} \frac{1}{(l_{1}l_{2})^{2}}\frac{1}{(l'_{1}l'_{2})^{2}} \frac{1}{(k k')^{4}} \max(\frac{1}{\log(T)^{2}}, \frac{1}{\log(T)^{2}T|kl_{1}-k'l'_{1}| }) 
\end{equation}
and
\begin{equation}
\label{chap4:eq118}
J_{2}(T) = \sum_{\substack{ \lceil \log(T) \rceil \geqslant k,k' \geqslant 1 \\ 1 \geqslant |kl_{1} - k'l'_{1}|\geqslant \frac{1}{T} }} \frac{1}{(l_{1}l_{2})^{2}}\frac{1}{(l'_{1}l'_{2})^{2}} \frac{1}{(k k')^{4}} \max(\frac{1}{\log(T)^{2}}, \frac{1}{\log(T)^{2}T|kl_{1}-k'l'_{1}| }) \textit{.}
\end{equation}
Yet, according to Lemma $\ref{chap4:lemme2}$, one has that: 
 \begin{equation}
\label{chap4:eq119}
J_{1}(T) \leqslant \frac{K}{T \log(T)^{2}} O(\log(T)^{2}) 
\end{equation}
and so $J_{1}(T) \rightarrow 0$ when $T \rightarrow \infty$. \\
Furthermore, with $P$ being a constant that can be chosen as large as one wants, one has, for $T$ large: 
\begin{equation}
\label{chap4:eq120}
J_{2}(T) \leqslant \frac{K}{T \log(T)^{2}} \int_{(l_{1},l_{2},l'_{1},l'_{2}) \in Dom(T)} \frac{1}{(l_{1}l_{2}l'_{1}l'_{2})^{2}} \frac{1}{l_{1} - l'_{1}} dl_{1} dl'_{1} dl_{2} dl'_{2} + o(1)  
\end{equation}
where 
\begin{align*}
Dom(T) = \{ (l_{1},l_{2},l'_{1},l'_{2}) \in \mathbb{R}^{4} \text{ } | \text{ } T \log(T) \geqslant l_{1} \geqslant P \text{, } \\
T \log(T) \geqslant l'_{1} \geqslant P-1 \text{, } T \log(T) \geqslant l'_{2},l_{2} > 0 \text{, } \\
l'_{1} +1 \geqslant l_{1} \geqslant l'_{1} + \frac{1}{T} \text{ and } l'_{1}l'_{2}\text{, } l_{1}l_{2} \geqslant C \}
\end{align*}
and $o(1)$ is a quantity that goes to $0$ when $T \rightarrow \infty$. \\
By integrating in $l_{2}$ and in $l'_{2}$, one gets from Equation ($\ref{chap4:eq120}$) that: 
\begin{equation}
\label{chap4:eq121}
J_{2}(T) \leqslant \frac{D}{T \log(T)^{2}} \int_{\substack{ T \log(T) \geqslant l_{1} \geqslant P \\ T \log(T) \geqslant l'_{1} \geqslant P-1 \\ l'_{1} +1 \geqslant l_{1}  \geqslant l'_{1} + \frac{1}{T}  }} \frac{1}{l_{1}l'_{1}} \frac{1}{l_{1}-l'_{1}}dl_{1}dl'_{1} + o(1) \textit{.}
\end{equation}
By using the fact that $0 < \frac{1}{l_{1}(l_{1}-l'_{1})} \leqslant \frac{1}{l'_{1}(l_{1}-l'_{1})}$ and by integrating on $l_{1}$, one gets with Equation $(\ref{chap4:eq121})$ that: 
\begin{equation}
\label{chap4:eq122}
J_{2}(T) \leqslant \frac{K}{T \log(T)^{2}} \int_{ T \geqslant l'_{1} \geqslant P-1} \frac{\log(T)}{(l'_{1})^{2}} dl'_{1} = O(\frac{1}{T \log(T)}) + o(1) \textit{.}
\end{equation}
So, we have finally that $J_{2}(T) \rightarrow 0$, when $T \rightarrow \infty$, and so does $J(T)$, $A_{2,2}(T)$ and $A_{2}(T)$.\\
By exchanging $l_{1}$ and $l_{2}$, we obtain the fact that $A_{1}(T) \rightarrow 0$ when $T \rightarrow \infty$. \\
Finally, with Equation ($\ref{chap4:eq108}$), one gets the wanted result. 
\end{proof}
%en fait, l'espérance est inutile, on peut passer directement à la norme 2 .... enlever cette partie sur l'espérance, juste le dire en note qu'on peut l'obtenir assez facilement ... 
%ensuite continuer sur le dernier truc: faire directement la norme 2 ... et il faudrait d'ailleurs enlever aussi la dépendance en t dans la somme ... 
\subsubsection{Proof of Proposition $\ref{chap4:prop11}$}
To prove Proposition $\ref{chap4:prop11}$, we are going to follow the same approach as was used just before. \\
Before entering into the proof of Proposition $\ref{chap4:prop11}$, we need the following preparatory lemma. Let us introduce: 
\begin{equation}
\label{chap4:eq123}
\tilde{G}_{4}(L,t,T) = \frac{1}{2 \pi^{4} \text{Covol}(L^{\perp})^{2}} \sum_{l \in J_{2}(L,T)} \frac{1}{(l_{1}l_{2})^{2}} \sum_{k=1}^{\infty} \frac{\cos (4 \pi k t l_{1})\cos (4 \pi k t l_{2}) }{k^{4}} \textit{. }
\end{equation}
Then we have the following lemma that in particular says that $\frac{G_{4}(L,t)}{V(L,t)}-\frac{\tilde{G}_{4}(L,t,T)}{V(L,t)}$ tends, in probability, towards $0$.
\begin{lemma}
\label{chap4:lemme11}
One has, when $T \rightarrow \infty$, 
$$\mathbb{E} \left( |\frac{\tilde{G}_{4}(L,t,T)}{V(L,t)}-\frac{G_{4}(L,t)}{V(L,t)}| \right) \rightarrow 0 \textit{.} $$
\end{lemma}
\begin{proof}
We use the same estimates as in the proof of Lemma $\ref{chap4:lemme10}$.
\end{proof}
We can now tackle the proof of Proposition $\ref{chap4:prop11}$.
\begin{proof}[Proof of Proposition $\ref{chap4:prop11}$]
According to Lemma $\ref{chap4:lemme11}$ and the Markov's inequality, we only need to prove that $\frac{\tilde{G}_{4}(L,t,T)}{V(L,t)}$ tends in probability towards $0$. \\
To do so, by using a usual trigonometric formula, we have that: 
\begin{equation}
\label{chap4:eq124}
\frac{\tilde{G}_{4}(L,t,T)}{V(L,t)} = D(U_{1}(L,t,T) + U_{2}(L,t,T)) 
\end{equation}
with $D > 0$ and where
\begin{equation}
\label{chap4:eq125}
U_{1}(L,t,T) =  \frac{1}{V(L,t)} \sum_{l \in J_{2}(L,T)} \frac{1}{(l_{1}l_{2})^{2}} \sum_{k=1}^{\infty} \frac{\cos (4 \pi k t (l_{1} - l_{2})) }{k^{4}}
\end{equation}
and
\begin{equation}
\label{chap4:eq126}
U_{2}(L,t,T) = \frac{1}{V(L,t)} \sum_{l \in J_{2}(L,T)} \frac{1}{(l_{1}l_{2})^{2}} \sum_{k=1}^{\infty} \frac{\cos (4 \pi k t (l_{1} + l_{2})) }{k^{4}} \textit{.}
\end{equation}
So, to prove that $\frac{\tilde{G}_{4}(L,t,T)}{V(L,t)}$ tends to $0$ in probability, we only need to show that the moments of order $2$ of $U_{1}(L,t,T)$ and $U_{2}(L,t,T)$ tend to $0$ in probability. We are going to prove this fact for $U_{1}(L,t,T)$ and the proof will be valid by symmetry for $U_{2}$ and we will so get the wanted result. \\
We have that: 
\begin{equation}
\label{chap4:eq127}
\mathbb{E}(U_{1}(L,t,T)^{2}) \leqslant D \sum_{l,l' \in J_{2}(L,T)} \frac{1}{(l_{1}l_{2}l'_{1}l'_{2})^{2}} \sum_{ k,k' \geqslant 1 } \frac{1}{(k k')^{4}} \min \left( \frac{1}{\log(T)^{2}}, \frac{h(k l,k l') +  h(k l,-k l')}{T \log(T)^{2}}  \right) 
\end{equation}
where $D > 0$ is a constant and 
\begin{equation}
\label{chap4:eq128}
h(l,l') = \frac{1}{|l_{1}-l_{2}-(l'_{1}-l'_{2})|} \textit{.}
\end{equation}
We are going to show that: 
\begin{equation}
\label{chap4:eq160}
Z(L,T) = \sum_{l,l' \in J_{2}(L,T)} \frac{1}{(l_{1}l_{2}l'_{1}l'_{2})^{2}} \sum_{ k,k' \geqslant 1 } \frac{1}{(k k')^{4}} \min \left( \frac{1}{\log(T)^{2}}, \frac{h(k l,k l') }{T \log(T)^{2}}  \right) \rightarrow 0
\end{equation}
when $T \rightarrow \infty$. Furthermore, the proof will still be valid if we exchange $l'$ and $-l'$. So we will get the wanted result due to Equation ($\ref{chap4:eq127}$). \\
One has that:
\begin{equation}
\label{chap4:eq129}
Z(L,T)  = \Delta_{1}(L,T) + \Delta_{2}(L,T) 
\end{equation}
where 
\begin{equation}
\label{chap4:eq130}
\Delta_{1}(L,T) = \sum_{l,l' \in J_{2}(L,T)} \frac{1}{(l_{1}l_{2}l'_{1}l'_{2})^{2}} \sum_{ \substack{ k,k' \geqslant 1 \\  k   \textit{ or }  k'   > \log(T) }} \frac{1}{(k k')^{4}} \min \left( \frac{1}{\log(T)^{2}}, \frac{h(k l,k l')}{T \log(T)^{2}} \right)
\end{equation}
and
\begin{equation}
\label{chap4:eq131}
\Delta_{2}(L,T) = \sum_{l,l' \in J_{2}(L,T)} \frac{1}{(l_{1}l_{2}l'_{1}l'_{2})^{2}} \sum_{ \substack{ k,k' \geqslant 1 \\  k   \textit{ and }  k'   \leqslant \log(T) }} \frac{1}{(k k')^{4}} \min \left( \frac{1}{\log(T)^{2}}, \frac{h(k l,k l') }{T \log(T)^{2}} \right) \textit{.}
\end{equation}
Yet, according to Lemma $\ref{chap4:lemme2}$, one has that: 
\begin{equation}
\label{chap4:eq132}
\Delta_{1}(L,T) \leqslant \frac{D}{\log(T)^{8}} (\int_{ \substack{ A \leqslant \lVert l \rVert \leqslant T  \\ |l_{1} l_{2} | \geqslant C }} \frac{1}{l_{1}^{2} l_{2}^{2}} dl_{1} dl_{2})^{2} + o(1) \underset{T \rightarrow \infty}{\rightarrow} 0
\end{equation}
where $D > 0$ is a constant. \\
So, we only need to prove that $\Delta_{2}(L,T) \rightarrow 0$ when $T \rightarrow \infty$. \\
Yet, one has also: 
\begin{equation}
\label{chap4:eq133}
\Delta_{2}(L,T) = \Delta_{2,1}(L,T) + \Delta_{2,2}(L,T) 
\end{equation}
where
\begin{equation}
\label{chap4:eq134}
\Delta_{2,1}(L,T)=  \sum_{l,l' \in J_{2}(L,T)} \frac{1}{(l_{1}l_{2}l'_{1}l'_{2})^{2}} \sum_{ \substack{ k,k' \geqslant 1 \\  k   \textit{ and }  k'   \leqslant \log(T)  \\h(k l,k l') \leqslant 1  }} \frac{1}{(k k')^{4}}  \frac{h(k l,k l') }{T \log(T)^{2}} \textit{ and }
\end{equation}
\begin{equation}
\label{chap4:eq135}
 \Delta_{2,2}(L,T)   = \sum_{l,l' \in J_{2}(L,T)} \frac{1}{(l_{1}l_{2}l'_{1}l'_{2})^{2}} \sum_{ \substack{ k,k' \geqslant 1 \\  k   \textit{ and }  k'   \leqslant \log(T)  \\  h(k l,k l') > 1  }} \frac{1}{(k k')^{4}}  \frac{h(k l,k l') }{T \log(T)^{2}} \textit{.}
\end{equation}
%\begin{equation}
%\label{chap4:eq136}
%\Delta_{2,3}(L,T) =  \Delta_{2,1}(L,T) =  \Delta_{2}(L,T) = \sum_{l,l' \in J_{2}(L,T)} \frac{1}{(l_{1}l_{2}l'_{1}l'_{2})^{2}} \sum_{ \substack{ k,k' \geqslant 1 \\  k   \textit{ and }  k'   \leqslant \log(T)  \\ h(k l , k l') > T }} \frac{1}{(k k')^{4}}  \frac{1}{\log(T)^{2}} \textit{.}
%\end{equation}
Yet, there exists $A > 0$ such that: 
\begin{equation}
\label{chap4:eq137}
\Delta_{2,1}(L,T) \leqslant \frac{K}{T \log(T)^{2}} \left( \int_{A \leqslant \lVert l \rVert \leqslant T \log(T)} \frac{1}{l_{1}^{2} l_{2}^{2}} \right)^{2} \rightarrow 0 
\end{equation}
when $T \rightarrow \infty$ and where $K$ is a positive constant. The right-hand side converges towards $0$ because of Lemma $\ref{chap4:lemme2}$. \\
Because of Equation $(\ref{chap4:eq133})$, it only remains to prove that $\Delta_{2,2}(L,T)$  converges towards $0$ when $T \rightarrow \infty$. \\
Yet, one has that: 
\begin{equation}
\label{chap4:eq138}
 \Delta_{2,2}(L,T) \leqslant \frac{K}{\log(T)^{2}} \int_{\substack{A \leqslant \lVert l \rVert, \lVert l' \rVert \leqslant T \log(T) \\ |l_{1}l_{2}|,|l'_{1}l'_{2}| \geqslant C \\ |l_{1}-l_{2}-(l'_{1}-l'_{2})| < 1 }} f(l,l')dl dl'
\end{equation}
where $K > 0$ and $f(l,l') = \left( \frac{1}{l_{1}l_{2}l'_{1}l'_{2}} \right)^{2}$. \\
Because of Equation ($\ref{chap4:eq138}$) and for symmetry reasons, it is enough to show that the following quantity converges towards $0$ when $T \rightarrow \infty$: 
\begin{equation}
\label{chap4:eq139}
J(T) = \frac{1}{\log(T)^{2}} \int_{(l,l') \in I(T) } f(l,l') dl dl'
\end{equation}
where 
\begin{align}
\label{chap4:eq140}
I(T)  = \{ (l,l') \in \mathbb{R}^{4} \textit{ } | & \textit{ } l_{1},l_{2},l'_{1},l'_{2} > 0 \textit{ , } \nonumber  \\
      & l_{1} > l_{2} \textit{ , } l'_{1} > l'_{2} \textit{ , } 1 > l_{1}-l_{2}-(l'_{1}-l'_{2}) > 0 \textit{ , } \nonumber \\
  &    A \leqslant  \lVert l \rVert, \lVert l' \rVert \leqslant T \log(T) \textit{ and } l_{1}l_{2},l'_{1}l'_{2} \geqslant C  \} \textit{.} 
\end{align}
Let us call $I_{1}(T)$ the set of all $(l,l')$ that belong to $I(T)$ and that verify $l_{1} - l_{2} \leqslant 2$. Let us call also $I_{2}(T) = I(T) - I_{1}(T)$. \\
Let us note that if $(l,l') \in I_{1}(T)$ then $l'_{1} - l'_{2} \leqslant 2$. \\
So, one has: 
\begin{equation}
\label{chap4:eq141}
J_{1}(T) = \frac{1}{\log(T)^{2}} \int_{(l,l') \in I_{1}(T) } f(l,l') dl dl' \underset{T \rightarrow \infty}{\rightarrow} 0 
\end{equation}
because $\int_{(l,l') \in I_{1}(T) } f(l,l') dl dl'$ is bounded ($l$ and $l'$ are close to the axis $y=x$). \\
Let us set: 
\begin{equation}
\label{chap4:eq142}
J_{2}(T)  =  \frac{1}{\log(T)^{2}} \int_{(l,l') \in I_{2}(T) } f(l,l') dl dl'
\end{equation}
so that 
\begin{equation}
\label{chap4:eq1105}
J(T) = J_{1}(T) + J_{2}(T) \textit{.}
\end{equation}
As a consequence from Equation ($\ref{chap4:eq141}$), it is enough to prove that $J_{2}(T) \rightarrow 0$ when $T \rightarrow \infty$ in order to prove that $\Delta_{2,2}(L,T)$ converges towards $0$ when $T \rightarrow \infty$. \\
It is easy to see that an $\textit{a priori}$ important part of the integral $\int_{(l,l') \in I_{2}(T) } f(l,l') dl dl'$ is the $l$ and $l'$ such that $l_{1}$ and $l'_{1}$ are large, for example larger than $\log(\log(T))$, and $l_{2}$ and $l'_{2}$ are small, for example smaller than $ \frac{1}{ \log(\log(T))}$. The rest of the integral, when divided by $\log(T)^{2}$, goes to $0$ when $T \rightarrow \infty$. \\
But we have also that $1 > l_{1}-l_{2}-(l'_{1}-l'_{2}) > 0$ because $(l,l') \in I(T)$. So, it implies in particular that, for $T$ large enough, $$ \frac{3}{2} > l_{1} - l'_{1} > - \frac{1}{2} \textit{.} $$ 
Hence, by integrating first in $l_{2}$ and in $l'_{2}$ and by using the fact that the lattice is admissible, we have that $J_{2}(T)$ is estimated as followed: 
\begin{equation}
\label{chap4:eq143}
J_{2}(T) \leqslant \frac{K_{1}}{\log(T)^{2}} \int_{\substack{ T \log(T) \geqslant l_{1},l'_{1} \geqslant \log(\log(T)) \\ \frac{3}{2} > l_{1} - l'_{1} > - \frac{1}{2}}} \frac{1}{l_{1} l'_{1} } dl_{1}dl'_{1} +o(1) 
\end{equation}
where $K_{1}$ is a positive constant (that depends on $C$). \\
From Equation ($\ref{chap4:eq143}$), by integrating, first, in $l_{1}$, second, in $l'_{1}$, one gets that: 
\begin{equation}
\label{chap4:eq144}
J_{2}(T) \leqslant \frac{K_{2}}{\log(T)}  +o_{T}(1) 
\end{equation}
where $K_{2} > 0$, which concludes the proof.

\end{proof}

\subsection{Conclusion}
We can now give the full proof of Theorem $\ref{chap4:thm1}$.
\begin{proof}[Proof of Theorem $\ref{chap4:thm1}$]
Proposition $\ref{chap4:prop2}$ gives us the first part of Theorem $\ref{chap4:thm1}$.
The second part of Theorem $\ref{chap4:thm1}$ is the consequence of Proposition $\ref{chap4:prop3}$ and of Proposition $\ref{chap4:prop8}$ (with $L$ being replaced by $L^{\perp}$ ; the validity of this last Proposition is a consequence of Proposition $\ref{chap4:prop9}$, Proposition $\ref{chap4:prop10}$ and Proposition $\ref{chap4:prop11}$).

\end{proof}

\section{Proof of Theorem 15}
The goal of this section is to establish Theorem $\ref{chap4:thm1000}$.
It is a natural extension of Theorem $\ref{chap4:thm1}$ and of a part of $\ref{chap4:thm2}$. We see that, in that case, the normalization is larger than before: at least $\log(t)$ whereas the normalization before was in $\sqrt{\log(t)}$. \\
To show this last theorem, we are going to proceed in three steps: first, we will show that the lower and upper estimates about $\tilde{V}(L,r)$ hold, second we will show that the lower and upper estimates about $\tilde{V}(L,r)$ hold and third we will conclude the proof of Theorem $\ref{chap4:thm1000}$ by using the third subsection of Section 4. \\
\subsection{Estimation of $\tilde{V}(L,r)$ in the typical case}
The goal of this subsection is to show the following proposition: 
\begin{prop}
\label{chap4:prop106}
For every $\epsilon > 0$, for a typical $L \in \mathscr{S}_{d}$, one has that 
\begin{equation}
\tilde{V}(L,r) = O(r^{2 d - 2 + \epsilon})
\end{equation}
and
\begin{equation}
\frac{\tilde{V}(L,r)}{r^{ d - 1}} \underset{r \rightarrow \infty}{\rightarrow} \infty \textit{.} 
\end{equation}
\end{prop}
The proof of this proposition relies heavily on $\cite{Skriganov}$ and we need to recall a fundamental theorem. 
\begin{defi}
A subgroup $G \subset SL_{d}(\mathbb{R})$ is called ergodic on the homogeneous space $\mathscr{S}_{d}$ if for every $G$-invariant measurable subset $A \subset \mathscr{S}_{d}$, $\mu_{d}(A) = 0$ or $\mu_{d}(A) = 1$ where $\mu_{d}$ is the unique Haar and probability measure on $\mathscr{S}_{d}$. 
\end{defi}
The Moore's ergodic theorem gives us in fact that $G$ is ergodic if, and only if, $G$ is not contained in any compact subgroup of $SL_{d}(\mathbb{R})$. As a consequence, $\Delta$ is ergodic and thus we have the following fundamental theorem
\begin{theorem}
\label{chap4:thm13}
Let $\psi$ a function integrable over $(\mathscr{S}_{d},\mu_{d})$. Then, for almost all $L \in \mathscr{S}_{d}$ (in the sense of the measure $\mu_{d}$), one has that 
$$\lim_{r \rightarrow \infty} \frac{1}{|\Delta_{r}|} \sum_{\delta \in \Delta_{r}} \psi(\delta L) = \int_{\mathscr{S}_{d}} \psi(L) d\mu_{d}(L) \textit{.} $$ 
\end{theorem}
We can now give the proof of Proposition $\ref{chap4:prop106}$.
\begin{proof}[Proposition $\ref{chap4:prop106}$]
First, one has that: 
\begin{equation}
\label{chap4:eq1034}
\tilde{V}(L,r) \leqslant \left( \sum_{\delta \in \Delta_{r}} \frac{1}{\lVert \delta L \rVert^{d}} \right)^{2} \textit{.}
\end{equation}
Yet, Lemma 3.2 from $\cite{Skriganov}$ gives us that for every $\epsilon > 0$, for a typical $L \in \mathscr{S}_{d}$, 
\begin{equation}
\label{chap4:eq1035}
\left( \sum_{\delta \in \Delta_{r}} \frac{1}{\lVert \delta L \rVert^{d}} \right) = O(r^{d-1+ \epsilon}) \textit{.}
\end{equation}
So, Equation ($\ref{chap4:eq1034}$) and Equation ($\ref{chap4:eq1035}$) give us the wanted first result of Proposition $\ref{chap4:prop106}$. \\
Second, one has that: 
\begin{equation}
\label{chap4:eq1036}
\sqrt{\frac{\tilde{V}(r,L)}{r^{ d - 1}}} \geqslant \frac{K_{1}}{r^{d-1}} \sum_{\delta \in \Delta_{r}} \frac{1}{\lVert \delta L \rVert^{d}} 
\end{equation}
where $K_{1} > 0$. We have obtained this last equation by using the concavity of the square root and because the cardinal number of $\Delta_{r}$ is of order $r^{d-1}$ (see Equation $(\ref{chap4:eq1300})$). \\
Let then $m \geqslant 1$. \\
One has that: 
\begin{equation}
\label{chap4:eq1037}
\frac{1}{r^{d-1}} \sum_{\delta \in \Delta_{r}} \frac{1}{\lVert \delta L \rVert^{d}}  \geqslant \frac{1}{r^{d-1}} \sum_{\delta \in \Delta_{r}} \min \left(m,\frac{1}{\lVert \delta L \rVert^{d}} \right) \textit{.} 
\end{equation}
Yet, $ L \longmapsto \min \left(m,\frac{1}{\lVert L \rVert^{d}} \right)$ is integrable over $\mathscr{S}_{d}$. %rappeler définiton de ergodic et tout ça ? 
So, one has for a typical $L \in \mathscr{S}_{d}$: 
\begin{equation}
\label{chap4:eq1038}
\frac{1}{r^{d-1}} \sum_{\delta \in \Delta_{r}} \min \left(m,\frac{1}{\lVert \delta L \rVert^{d}} \right) \underset{r \rightarrow \infty}{\rightarrow} K_{2} \int_{\mathscr{S}_{d}} \min(m, \frac{1}{\lVert L \rVert^{d}}) d\mu_{d}  
\end{equation}
where $K_{2} > 0$ does not depend on $m$. \\
By using Equation $(\ref{chap4:eq1036})$, Equation $(\ref{chap4:eq1037})$ and Equation $(\ref{chap4:eq1038})$, one gets that for every $m \geqslant 1$, for a typical $L  \in \mathscr{S}_{d}$:
\begin{equation}
\label{chap4:eq1039}
\liminf_{r \rightarrow \infty} \sqrt{\frac{\tilde{V}(r,L)}{r^{ d - 1}}} \geqslant K \int_{\mathscr{S}_{d}} \min(m, \frac{1}{\lVert L \rVert^{d}}) d\mu_{d} 
\end{equation}
where $K >0$ and $m > 0$. \\
By making $m \rightarrow \infty$, by using Fatou's lemma and by using the fact that $$\int_{\mathscr{S}_{d}}\frac{1}{\lVert L \rVert^{d}} d\mu_{d} = \infty \textit{, }$$ one gets the wanted result.
\end{proof}
\subsection{Result about $V(L,t)$ in the typical case}
The goal of this subsection is to prove the following proposition: 
\begin{prop}
\label{chap4:prop107}
For every $\epsilon > 0$, for $L$ a typical lattice, one has that
\begin{equation*}
V(L,t) = O(\log(t)^{2+\epsilon})
\end{equation*}
and
\begin{equation*}
\frac{V(L,t)}{\log(t)^{2}} \underset{t \rightarrow \infty}{\rightarrow} \infty \textit{.}
\end{equation*}
\end{prop}
In fact, the most important part of $V(L,t)$ are the terms $l_{1}^{2} l_{2}^{2}$ that are the smallest possible. So, we need to know more about how small can be $|l_{1} l_{2}|$. In fact, we have the following result: 
\begin{theorem}[\cite{sprindzhuk1979metric}, \cite{Skriganov}]
\label{chap4:thm12}
For a typical $L \in \mathscr{S}_{2}$, there exists a sequence $(l_{n})_{n \in \mathbb{N}}$ such that $$ \rVert l_{n} \lVert \underset{n \rightarrow \infty}{\rightarrow} \infty \textit{ and } \log(\lVert l_{n} \rVert) |(l_{n})_{1} (l_{n})_{2}| \underset{n \rightarrow \infty}{\rightarrow} 0 \textit{.} $$ 
Furthermore, for all $\alpha > 0$, for a typical $L \in \mathscr{S}_{2}$, there exists $C > 0$ such that for all $l \in L-\{0 \}$,
$$ |l_{1}l_{2}| \geqslant C | \log(\lVert l \rVert)|^{-1 - \alpha} \textit{.}$$
\end{theorem}
As a consequence of Theorem $\ref{chap4:thm12}$, we see that, to establish Proposition $\ref{chap4:prop107}$, the following lemma will be convenient: 
\begin{lemma}
\label{chap4:lemme13}
For all $C > 0$, for all $A > 0$, for all $\alpha > 0$, one has that: 
$$\int_{\substack{ l \in \mathbb{R}^{2} \\ A \leqslant \lVert l \rVert \leqslant t  \\ |\text{Num}(l)| \geqslant C | \log(\lVert l \rVert) |^{-1-\alpha}}} \frac{1}{l_{1}^{2}  l_{2}^{2}} dl_{1} dl_{2} = O (\log(t)^{2+ \alpha}) \textit{.}$$
\end{lemma}
The proof of this lemma is basically the same as the proof of Lemma $\ref{chap4:lemme2}$. 
\begin{proof}
First, let us say that it is enough to prove the result for $A$ fixed and large enough. Then, let us set: 
\begin{equation}
\label{chap4:eq1040}
J(t) = \int_{\substack{ l \in \mathbb{R}^{2} \\ A \leqslant \lVert l \rVert \leqslant t  \\ |\text{Num}(l)| \geqslant C | \log(\lVert l \rVert) |^{-1-\alpha}}} \frac{1}{l_{1}^{2}  l_{2}^{2}} dl_{1} dl_{2} \textit{.}
\end{equation}

By passing into polar coordinates $(r,\theta)$ and by using the symmetries, one has: 
\begin{equation}
\label{chap4:eq1041}
J(t) = 8 \int_{\substack{A \leqslant r \leqslant t \\ 0 \leqslant \theta \leqslant \frac{\pi}{4} \\ \sin(2 \theta) \geqslant \frac{2 C \log(r)^{-1-\alpha} }{r^{2}} }} \frac{1}{r^{3}} \frac{4}{\sin(2 \theta)^{2}} dr d\theta \textit{.}
\end{equation} %16
By making the changes of variable $\theta' = 2 \theta$ and, then, $u = \tan(\theta')$, one gets from Equation $(\ref{chap4:eq11})$ and by taking $A$ large enough: 
\begin{equation}
\label{chap4:eq1042}
J(t) = 16 \int_{A \leqslant r \leqslant t^{1+ \epsilon_{1}}} \frac{1}{r^{3}} (-1 + \frac{r^{2} \log(r)^{1+\alpha}}{2 C} \sqrt{1 - (\frac{2 C}{\log(r)^{1+\alpha} r^{2}})^{2}}) dr \textit{.} 
\end{equation}
From Equation ($\ref{chap4:eq1042}$), for all $A$ large enough, one has that: 
\begin{equation}
\label{chap4:eq1043}
J(t) = O (\log(t)^{2+\alpha}) 
\end{equation}
when $t \rightarrow \infty$. So we get the wanted result.
\end{proof}
We can now prove Proposition $\ref{chap4:prop107}$.
\begin{proof}[Proof of Proposition $\ref{chap4:prop107}$]
Let $\epsilon > 0$. For a typical $L \in \mathscr{S}_{2}$, there exists $C > 0 $ such that for all $l \in L-\{0 \}$
\begin{equation}
\label{chap4:eq1044}
|l_{1}l_{2}| \geqslant C | \log(\lVert l \rVert)|^{-1 - \epsilon}
\end{equation}
and there exists a sequence $(l_{n})_{n \in \mathbb{N}} \in L^{\mathbb{N}}$ such that 
\begin{equation}
\label{chap4:eq1045}
 \rVert l_{n} \lVert \underset{n \rightarrow \infty}{\rightarrow} \infty \textit{ and }
 \end{equation}
 \begin{equation}
 \label{chap4:eq1046} \log(\lVert l_{n} \rVert) |(l_{n})_{1} (l_{n})_{2}| \underset{n \rightarrow \infty}{\rightarrow} 0 \textit{.} 
 \end{equation}
So, first, one has that there exist $A,D > 0$ such that: 
\begin{equation}
\label{chap4:eq1047}
V(L,t) \leqslant D \int_{\substack{ l \in \mathbb{R}^{2} \\ A \leqslant \lVert l \rVert \leqslant t  \\ |\text{Num}(l)| \geqslant C | \log(\lVert l \rVert) |^{-1-\epsilon}}} \frac{1}{l_{1}^{2}  l_{2}^{2}} dl_{1} dl_{2} = O (\log(t)^{2+ \epsilon})
\end{equation}
according to the definition of $V(L,t)$ and because of Equation $(\ref{chap4:eq1044})$. \\
So, because of Lemma $\ref{chap4:lemme13}$, one gets that: 
\begin{equation}
V(L,t) = O (\log(t)^{2+ \epsilon}) \textit{.}
\end{equation}
Furthermore, a consequence of Equation $(\ref{chap4:eq1045})$ and of Equation $(\ref{chap4:eq1046})$ is the fact that: 
\begin{equation}
\label{chap4:eq1048}
\liminf_{t \rightarrow \infty} \frac{V(L,t)}{\log(t)^{2}} = \infty
\end{equation}
also due to the definition of $V(L,t)$. 
\end{proof}
We can now conclude the proof of Theorem $\ref{chap4:thm1000}$. 
\subsection{Conclusion}
\begin{proof}[Proof of Theorem $\ref{chap4:thm1000}$]
The first assertion of Theorem $\ref{chap4:thm1000}$ is proven by Proposition $\ref{chap4:prop106}$.\\
The second assertion of Theorem $\ref{chap4:thm1000}$ is proven by Proposition $\ref{chap4:prop107}$.\\
The third part of Theorem $\ref{chap4:thm1000}$, that concerns the convergence in distribution, is shown as the second part of Theorem $\ref{chap4:thm1}$ (see 4.2 and 4.3). 
\end{proof}

\section{Proof of Theorem 16}
In this section, we are going to show Theorem $\ref{chap4:thm3}$. Let $x \in \mathbb{R}$ and $a > 0$. \\
Instead of considering $t$, we can consider $ \frac{t}{a}$. So, we can suppose, and we are going to make this assumption in the rest of this section, that $a=1$ in the study of $\frac{\mathcal{R}(t \text{Rect}(a,a) + (x,x) ,\mathbb{Z}^{2})}{t}$. \\
Now we are going to give a simple expression of $\frac{\mathcal{R}(t \text{Rect}(1,1) + (x,x) ,\mathbb{Z}^{2})}{t}$.
%donner plan de la preuve
\subsection{An expression of $\frac{\mathcal{R}(t \text{Rect}(1,1) + (x,x) ,\mathbb{Z}^{2})}{t}$}
The main objective of this subsection is to prove the following proposition: 
\begin{prop}
\label{chap4:prop100}
We have for every $x \in \mathbb{R}^{2}$, for every $t > 0$, that 
\begin{equation}
\label{chap4:eq1001}
\frac{\mathcal{R}(t \text{Rect}(1,1) + (x,x) ,\mathbb{Z}^{2})}{t} = \frac{(\lfloor t + x \rfloor - \lceil -t + x \rceil +1)^{2} - 4 t^{2}}{t} \textit{.}
\end{equation}
\end{prop}
The proof is quite straightforward.
\begin{proof}
Let $x \in \mathbb{R}^{2}$. Let $t > 0$. \\
One has that: 
\begin{align}
\label{chap4:eq1002}
N(t \text{Rect}(1,1) + (x,x) ,\mathbb{Z}^{2})& = \left( \sum_{ \substack{ (n_{1},n_{2}) \in \mathbb{Z}^{2} \\ -t+x \leqslant n_{1}  \leqslant t+x \\ -t+x \leqslant n_{2}  \leqslant t+x }} 1 \right) \nonumber \\
& = (\lfloor t + x \rfloor - \lceil -t + x \rceil +1)^{2} \textit{.}
\end{align}
Furthermore, one has that: 
\begin{equation}
\label{chap4:eq1003}
\text{Area}(t \text{Rect}(1,1))= 4 t^{2} \textit{.}
\end{equation}
So, Equation $(\ref{chap4:eq1002})$ and Equation $(\ref{chap4:eq1003})$ and the definition of $\mathcal{R}(t \text{Rect}(1,1) + (x,x) ,\mathbb{Z}^{2})$ give us Equation $(\ref{chap4:eq1001})$.
\end{proof}
With Equation $(\ref{chap4:eq1001})$, one has that: 
\begin{equation}
\label{chap4:eq1004}
\frac{\mathcal{R}(t \text{Rect}(1,1) + (x,x) ,\mathbb{Z}^{2})}{t} = (\lfloor t + x \rfloor - \lceil -t + x \rceil +1 - 2 t )\frac{(\lfloor t + x \rfloor - \lceil -t + x \rceil +1 + 2 t )}{t} \textit{.}
\end{equation}
Thanks to this last remark, the asymptotical study of $\frac{\mathcal{R}(t \text{Rect}(1,1) + (x,x) ,\mathbb{Z}^{2})}{t}$ with $t$ distributed on $[0,T]$ according to the probability measure $\frac{1}{T} \rho(\frac{t}{T}) dt$ is going to be reduced to the study of a simpler quantity. It is the object of the next subsection.
\subsection{Reduction of the study of $\frac{\mathcal{R}(t \text{Rect}(1,1) + (x,x) ,\mathbb{Z}^{2})}{t}$}
The main objective of this subsection is to prove the following proposition: 
\begin{prop}
\label{chap4:prop101}
For every $g \in C_{c}(\mathbb{R})$, 
\begin{equation}
\label{chap4:eq1005}
\int_{t=0}^{T} (g\left( \frac{\mathcal{R}(t \text{Rect}(1,1) + (x,x) ,\mathbb{Z}^{2})}{t} \right) - g\left( \Delta(t,x) \right)) \frac{1}{T} \rho(\frac{t}{T}) dt \underset{T \rightarrow \infty}{\rightarrow} 0
\end{equation}
where $\Delta(t,x)$ is defined by 
\begin{equation}
\label{chap4:eq1111}
\Delta(t,x) = 4 (\lfloor t + x \rfloor - \lceil -t + x \rceil +1 - 2 t ) \textit{.}
\end{equation}
\end{prop}
The proof is quite straightforward and lie on the definitions of $\lfloor \cdot \rfloor$ and of $\lceil \cdot \rceil$.
\begin{proof}
For every $t > 0$, for every $x \in \mathbb{R}$, 
\begin{equation}
\label{chap4:eq1006}
t+x - 1< \lfloor t + x \rfloor \leqslant t+x
\end{equation}
and 
\begin{equation}
\label{chap4:eq1007}
-t+x \leqslant \lceil -t + x \rceil < -t+x +1 \textit{.}
\end{equation}
From Equation $(\ref{chap4:eq1006})$ and Equation $(\ref{chap4:eq1007})$, one gets that: 
\begin{equation}
\label{chap4:eq1008}
\frac{4 t - 1}{t}< \frac{(\lfloor t + x \rfloor - \lceil -t + x \rceil +1 + 2 t )}{t} \leqslant \frac{4 t + 1}{t}
\end{equation}
So, from this last equation, one has that, when $t \rightarrow \infty$, 
\begin{equation}
\label{chap4:eq1009}
\frac{(\lfloor t + x \rfloor - \lceil -t + x \rceil +1 + 2 t )}{t} \rightarrow 4 \textit{.}
\end{equation}
So, one has that: 
\begin{equation}
\label{chap4:eq1010}
\frac{\mathcal{R}(t \text{Rect}(1,1) + (x,x) ,\mathbb{Z}^{2})}{t} - \Delta(t,x) \rightarrow 0
\end{equation}
when $t \rightarrow \infty$ because $(\lfloor t + x \rfloor - \lceil -t + x \rceil +1 - 2 t )$ is bounded. \\
Now, the end of this proof is quite straightforward. Indeed, one has, for every $0 < \kappa < \frac{1}{2}$: 
\begin{align}
\label{chap4:eq1011}
&| \int_{t=0}^{T} (g\left( \frac{\mathcal{R}(t \text{Rect}(1,1) + (x,x) ,\mathbb{Z}^{2})}{t} \right) - g\left( \Delta(t,x) \right)) \frac{1}{T} \rho(\frac{t}{T}) dt |  \leqslant   2 \lVert  g \rVert_{\infty} \int_{0}^{\kappa} \rho(t) dt  \nonumber \\
& + \int_{\kappa T}^{T} |g\left( \frac{\mathcal{R}(t \text{Rect}(1,1) + (x,x) ,\mathbb{Z}^{2})}{t} \right) - g\left( \Delta(t,x) \right)| \frac{1}{T} \rho(\frac{t}{T}) dt
\end{align}
and, because $g \in C_{c}(\mathbb{R})$, it is a uniformly continuous function and one has, from Equation $(\ref{chap4:eq1011})$, that
\begin{equation}
\label{chap4:eq1012}
\limsup_{T \rightarrow \infty} \int_{\kappa T}^{T} |g\left( \frac{\mathcal{R}(t \text{Rect}(1,1) + (x,x) ,\mathbb{Z}^{2})}{t} \right) - g\left( \Delta(t,x) \right)| \frac{1}{T} \rho(\frac{t}{T}) dt = 0
\end{equation}
because, also, of Equation $(\ref{chap4:eq1010})$. \\
So, Equation $(\ref{chap4:eq1011})$ and Equation $(\ref{chap4:eq1012})$ give us that for every $0 < \kappa \leqslant \frac{1}{2}$: 
\begin{equation}
\label{chap4:eq1013}
\limsup_{T \rightarrow \infty} | \int_{t=0}^{T} (g\left( \frac{\mathcal{R}(t \text{Rect}(1,1) + (x,x) ,\mathbb{Z}^{2})}{t} \right) - g\left( \Delta(t,x) \right)) \frac{1}{T} \rho(\frac{t}{T}) dt |  \leqslant   2 \lVert  g \rVert_{\infty} \int_{0}^{\kappa} \rho(t) dt   \textit{.}
\end{equation}
By making $\kappa$ go to $0$, one gets the wanted result.
\end{proof}
Proposition $\ref{chap4:prop101}$ enables us to reduce the asymptotic study of $\frac{\mathcal{R}(t \text{Rect}(1,1) + (x,x) ,\mathbb{Z}^{2})}{t}$ to the asymptotic study of $\Delta(t,x)$. In the next subsection, we are going to show that $\Delta(t,x)$ converges in distribution and exhibit the limit distribution and its moments.
\subsection{Convergence in distribution of $\Delta(t,x)$}
The goal of this subsection is to prove the following proposition: 
\begin{prop}
\label{chap4:prop105}
For all $x \in \mathbb{R}$, when $t \in [0,T]$ is distributed according to the probability density $\frac{1}{T} \rho(T \cdot)$ on $[0,T]$ then, when $T \rightarrow \infty$, $\Delta(t,x)$ converges in distribution. Furthermore, the limit distribution $\beta$ has a compact support included in $[-2,4]$ and for every $k \in \mathbb{N}$, one has that 
$$\int_{x \in \mathbb{R}} x^{k} d \beta(x) = a_{k} $$ 
where
$$
a_{k} = \frac{4^{k}(1 + (-1)^{k})(y^{k+1} + (1-y)^{k+1})}{2(k+1)} 
$$
with $y = |t_{2,0} - t_{1,0}|$ where $t_{2,0}$ is the first $t \geqslant 0$ such that $ -t + x \in \mathbb{Z}$ and $t_{1,0}$ is the first $t \geqslant 0$ such that $ t + x \in \mathbb{Z}$.
\end{prop}
We are going to show this proposition in three steps. The first time, and the next subsubsection is dedicated to it, consists in calculating the limit, when $T \rightarrow \infty$, of all entire moments of $\Delta(t,x)$ when $\rho = \mathbf{1}_{[0,1]}$. The second time consists in showing that these limits define a unique probability distribution over $\mathbb{R}$. The last subsection is dedicated to the conclusion of the proof. So, basically, we are going to use the method of moments to show Proposition $\ref{chap4:prop105}$.
\subsubsection{Calculation of limits of moments of $\Delta(t,x)$}
Before stating the main proposition of this section, we need to make some observations and put in place some notations.\\
Let us call $t_{1,0} < \cdots < t_{1,l}$ the different steps $t \in [0,T]$ such that $t + x \in \mathbb{Z}$. \\
In the same way, let us call $t_{2,0} < \cdots < t_{2,h}$ the different steps $t \in [0,T]$ such that $-t + x \in \mathbb{Z}$. \\
Let us observe that for every $i \in \{0, \cdots, l-1 \}$, $t_{1,i+1} - t_{1,i} = 1$ and that for every $j \in \{0, \cdots, h-1 \}$, $t_{2,j+1} - t_{2,j} = 1$. \\
As a consequence, one has necessarily that $t_{2,0} \in [t_{1,0},t_{1,1}[$ or $t_{1,0} \in [t_{2,0},t_{2,1}[$ and $h=l$ or $h=l-1$ or $h=l+1$. \\
Let us set: 
\begin{equation}
\label{chap4:eq1014}
y = |t_{1,0} - t_{2,0}| \textit{.}
\end{equation}
Then, the main proposition of this section is the following proposition: 
\begin{prop}
\label{chap4:prop102}
For every $k \in \mathbb{N}$, when $\rho = \mathbf{1}_ {[0,1]}$, one has that: 
\begin{equation}
\label{chap4:eq1015}
\lim_{T \rightarrow \infty} \left( \mathbb{E}((\Delta(t,x))^{k}) \right) = a_{k}
\end{equation}
where 
\begin{equation}
\label{chap4:eq1023}
a_{k} = \frac{4^{k}(1 + (-1)^{k})(y^{k+1} + (1-y)^{k+1})}{2(k+1)} \textit{.}
\end{equation}
\end{prop}
The proof consists basically in cutting the interval $[0,T]$ into subintervals where all the quantities that intervene in the calculus can be expressed simply. 
\begin{proof}
Let $k \geqslant 0$ and let us suppose that $\rho = \mathbf{1}_{[0,1]}$. \\
By symmetry, we can, and we will, also suppose that $t_{2,0} \in [t_{1,0},t_{1,1}[$.  \\
One has that:
\begin{equation}
\label{chap4:eq1110}
-4 \leqslant \Delta(t,x) \leqslant 4
\end{equation}  
for all $t \in \mathbb{R}$ and $x \in \mathbb{R}$.\\
Consequently, we can suppose that  
\begin{equation}
\label{chap4:eq1016}
\mathbb{E}((\Delta(t,x))^{k}) = \sum_{i=0}^{h-1} \int_{t_{1,i}}^{t_{2,i}} \Delta(t,x)^{k} \frac{1}{T} dt + \int_{t_{2,i}}^{t_{1,i+1}} \Delta(t,x)^{k} \frac{1}{T} dt 
\end{equation}
even if it means neglecting the rest of the integral that is calculated on a union of two intervals of respective lengths at most $2$ and so the corresponding term, because of Equation $(\ref{chap4:eq1110})$, is a $O(\frac{1}{T})$. \\
Let $i \in \{0, \cdots, h-1\}$. \\
Then one has: 
\begin{equation}
\label{chap4:eq1017}
\int_{t_{1,i}}^{t_{2,i}} \Delta(t,x)^{k} \frac{1}{T} dt = \int_{t_{1,i}}^{t_{2,i}} 4^{k} (t_{1,i} + t_{2,i}  -2 t)^{k} \frac{1}{T} dt
\end{equation}
according to the Equation $(\ref{chap4:eq1111})$. \\
So, one gets that: 
\begin{equation}
\label{chap4:eq1018}
\int_{t_{1,i}}^{t_{2,i}} \Delta(t,x)^{k} \frac{1}{T} dt = \frac{4^{k}}{2 T(k+1)}( y^{k+1} - (-y)^{k+1} ) 
\end{equation}
because for all $i \in \{0,\cdots,h-1 \}$, $ y= t_{2,0}-t_{1,0} = t_{2,i} - t_{1,i} $. \\
So, one gets that: 
\begin{equation}
\label{chap4:eq1019}
\int_{t_{1,i}}^{t_{2,i}} \Delta(t,x)^{k} \frac{1}{T} dt = \frac{4^{k}}{2 T(k+1)} y^{k+1} (1 + (-1)^{k})  \textit{.}
\end{equation}
In a similar way, one gets that: 
\begin{equation}
\label{chap4:eq1020}
\int_{t_{2,i}}^{t_{1,i+1}} \Delta(t,x)^{k} \frac{1}{T} dt =  \frac{4^{k}}{2 T(k+1)} (1-y)^{k+1} (1 + (-1)^{k})  
\end{equation}
because $1-y =  t_{1,i+1} - t_{2,i}$. \\
So, with Equation $(\ref{chap4:eq1016})$, Equation $(\ref{chap4:eq1019})$ and Equation $(\ref{chap4:eq1020})$, one gets that: 
\begin{equation}
\label{chap4:eq1021}
\mathbb{E}((\Delta(t,x))^{k}) = \sum_{i=0}^{h-1} \frac{4^{k}}{2 T(k+1)}(y^{k+1}+(1-y)^{k+1}) (1 + (-1)^{k})  \textit{.}
\end{equation}
By using the fact that $ \lim_{T \rightarrow \infty} \frac{h}{T} = 1$, one gets from Equation $(\ref{chap4:eq1021})$ that: 
\begin{equation}
\label{chap4:eq1022}
\mathbb{E}((\Delta(t,x))^{k}) =  \frac{4^{k} h}{2 T(k+1)}(y^{k+1}+(1-y)^{k+1}) (1 + (-1)^{k}) \underset{T \rightarrow \infty}{\rightarrow} \frac{4^{k} }{2(k+1)}(y^{k+1}+(1-y)^{k+1}) (1 + (-1)^{k})  \textit{.}
\end{equation}
%Dire par symétrie on fait le cas $t_{1,0} \leqslant t_{2,0}$. = on trouve le même résultat, on fait juste inégalité stricte. Puis expliquer que on se stoppe à t_{1,h} machin, le reste tend vers 0 de toute manière. Puis faire le calcul.
\end{proof}
In the next subsubsection, we are going to see that there exists a unique probability distribution over $\mathbb{R}$ such that the entire moments are given by the $a_{k}$.
\subsubsection{Existence and uniqueness of the distribution whose moments are given by the $a_{k}$}
The main objective of this subsubsection is to prove the following proposition: 
\begin{prop}
\label{chap4:prop103}
There exists a unique probability distribution $\beta$ over $\mathbb{R}$ such that for all $k \in \mathbb{N}$, 
$$\int_{x \in \mathbb{R}} x^{k} d\beta(x) = a_{k} \textit{.} $$
\end{prop}
To prove this proposition, we need to recall the following theorem.
\begin{theorem}
\label{chap4:thm7}
Let $(\alpha_{k})_{k \in \mathbb{N}}$ be a sequence of real numbers such that the power series $\sum_{k \geqslant 0} \frac{\alpha_{k}}{k!} z^{k}$ has a positive radius of convergence. Then, there exists at most one probability measure $\beta$ over $\mathbb{R}$ such that for all $ k\in \mathbb{N}$, 
$$ \alpha_{k} =  \int_{x \in \mathbb{R}} x^{k} d\beta(x) \textit{.}$$
\end{theorem}
We can now prove the proposition $\ref{chap4:prop103}$.
\begin{proof}[Proof of Proposition $\ref{chap4:prop103}$]
One has:
\begin{equation}
\label{chap4:eq1024}
\frac{\frac{a_{k+1}}{(k+1)!}}{\frac{a_{k}}{k!}} \underset{k \rightarrow \infty} { \rightarrow} 0
\end{equation}
according to Equation ($\ref{chap4:eq1023}$). \\
So, the ratio test gives us that $\sum_{k \geqslant 0} \frac{a_{k}}{k!} z^{k}$ has a radius of convergence that is infinite. \\
As a consequence, Theorem $\ref{chap4:thm7}$ applies here and gives us that there is at most one probability measure $\beta$ over $\mathbb{R}$ whose entire moments are given by the $(a_{k})_{k \in \mathbb{N}}$. \\
Furthermore, $$-4 \leqslant \Delta(t,x) \leqslant 4$$ for all $t \geqslant 0$ and for all $x \in \mathbb{R}$.\\
So, we have that for every $x \in \mathbb{R}$, $(\mathbb{P}_{\Delta(\cdot,x)})_{T > 0}$ is tight (recall that the distribution of $t$ depends on $T$), which means here that for every $\epsilon > 0$, there exists a compact set $W_{\epsilon}$ of $\mathbb{R}$ such that for every $T >0$, $$\mathbb{P}_{\Delta(\cdot,x)}\big(W_{\epsilon} \big) \geqslant 1 - \epsilon \textit{.}$$ 
As a consequence, Prokhorov's theorem gives us the existence of the probability measure $\beta$ whose entire moments are given by the $(a_{k})_{k \in \mathbb{N}}$. 
\end{proof}
\subsubsection{Conclusion of the proof of Proposition $\ref{chap4:prop105}$}
To conclude the proof of Proposition $\ref{chap4:prop105}$, we need to recall one theorem and to prove one lemma.
\begin{theorem}
\label{chap4:thm8}
Let $X$ be a real random variable characterized by its entire moments and $(X_{n})_{n \in \mathbb{N}}$ be a sequence of real random variables such that for all $k \in \mathbb{N}$, 
$$E(X_{n}^{k}) \underset{n \rightarrow \infty}{\rightarrow} E(X^{k}) \textit{.}$$
Then the sequence $(X_{n})$ converges in distribution towards $X$. 
\end{theorem}
The following lemma is in fact taken from $\cite{bleher1992distribution}$ (but it was not formulated as it, see proof of theorem 4.3). We are going to give the proof of this lemma here for completeness. 
\begin{lemma}
\label{chap4:lemme12}
Let $F$ be a real measurable function from $\mathbb{R}_{+}$. \\
Assume that there exists a probability measure $\mu$ over $\mathbb{R}$ such that for every $g \in C_{b}(\mathbb{R})$,
$$\lim_{T \rightarrow \infty} \frac{1}{T} \int_{0}^{T} g(F(t)) dt = \int_{x \in \mathbb{R}} g(x) d\mu(x) \textit{.}$$
Then, for every probability density $\rho$ over $[0,1]$, for every $g \in C_{b}(\mathbb{R})$, one has
\begin{equation}
\label{chap4:eq1025}
\lim_{T \rightarrow \infty} \frac{1}{T} \int_{0}^{T} g(F(t)) \rho(\frac{t}{T}) dt = \int_{x \in \mathbb{R}} g(x) d\mu(x) \textit{.}
\end{equation}
\end{lemma}
\begin{proof}
Assume first that $\rho$ is a step-wise function consisting of a finite number of steps. By linearity of the Equation $(\ref{chap4:eq1025})$, it is enough, in this case, to prove ($\ref{chap4:eq1025}$) for the function 
$$ \rho(x) = \frac{1}{b-a} \mathbf{1}_{[a,b]}(x) $$ where $0 \leqslant a < b \leqslant 1$, $\textit{id est}$ a one-step function. \\
Let $g \in C_{b}(\mathbb{R})$. \\
Because of the assumption of Lemma $\ref{chap4:lemme12}$, one has: 
\begin{align}
\label{chap4:eq1026}
\frac{1}{T} \int_{0}^{T} g(F(t)) \rho(\frac{t}{T}) dt & = \frac{1}{T(b-a)} \int_{aT}^{bT} g(F(t)) dt \nonumber \\ 
 & = \frac{1}{b-a}( b \frac{1}{b T} \int_{0}^{bT} g(F(t)) dt -  a \frac{1}{a T} \int_{0}^{aT} g(F(t)) dt ) \nonumber \\
 & \underset{T \rightarrow \infty}{ \rightarrow} \frac{1}{b-a} (b \int_{x \in \mathbb{R}} g(x) d\mu(x) - a \int_{x \in \mathbb{R}} g(x) d\mu(x)) = \int_{x \in \mathbb{R}} g(x) d\mu(x)
\end{align}
So we have Equation $(\ref{chap4:eq1025})$ in this case. \\
The general case follows now by using the fact that for every probability density $\rho$ over $[0,1]$, for every $\epsilon > 0$, there exists $\rho_{\epsilon}$, a step-wise function consisting of a finite number of steps, such that 
$$\int_{x \in [0,1]} |\rho(x) - \rho_{\epsilon}(x)| dx \leqslant \epsilon \textit{.} $$
\end{proof}
We can now prove Proposition $\ref{chap4:prop105}$. \\
\begin{proof}[Proof of Proposition $\ref{chap4:prop105}$]
Thanks to Lemma $\ref{chap4:lemme12}$, we only need to prove Proposition $\ref{chap4:prop105}$ in the case where $\rho = \mathbf{1}_{[0,1]}$. \\
Furthermore, Proposition $\ref{chap4:prop102}$, Proposition $\ref{chap4:prop103}$ and Theorem $\ref{chap4:thm8}$ gives us that that $\Delta(\cdot,x)$ converges in distribution when $T \rightarrow \infty$ and the moments of the limit distribution is given by the $a_{k}$. \\
Finally, we recall that $-4 \leqslant \Delta(t,x) \leqslant 4$ for all $t \geqslant 0$, for all $x \in \mathbb{R}$ and so $\beta$ has its support compact and included in $[-4,4]$. 
\end{proof}
\subsection{Conclusion}
We can now prove Theorem $\ref{chap4:thm3}$.
\begin{proof}[Proof of Theorem $\ref{chap4:thm3}$]
Thanks to Lemma $\ref{chap4:lemme12}$, we only need to prove Theorem $\ref{chap4:thm3}$ in the case where $\rho = \mathbf{1}_{[0,1]}$ and $a=1$. \\
Proposition $\ref{chap4:prop101}$ gives us that, if $\Delta(\cdot,x)$ converges in distribution, then it is also the case of $\frac{\mathcal{R}(t \text{Rect}(1,1) + (x,x), \mathbb{Z}^{2})}{t}$ and the limit distribution is the same. \\
Finally, Proposition $\ref{chap4:prop105}$ gives the wanted result. 
\end{proof}
\bibliographystyle{plain}
\bibliography{bibliographie}

\begin{thebibliography}{10}

\bibitem{beck2010randomness}
J{\'o}zsef Beck.
\newblock Randomness of the square root of 2 and the giant leap, part 1.
\newblock {\em Periodica Mathematica Hungarica}, 60(2):137--242, 2010.

\bibitem{bleher1992distribution}
Pavel Bleher et~al.
\newblock On the distribution of the number of lattice points inside a family
  of convex ovals.
\newblock {\em Duke Mathematical Journal}, 67(3):461--481, 1992.

\bibitem{bleher1993distribution}
Pavel~M Bleher, Zheming Cheng, Freeman~J Dyson, and Joel~L Lebowitz.
\newblock Distribution of the error term for the number of lattice points
  inside a shifted circle.
\newblock {\em Communications in mathematical physics}, 154(3):433--469, 1993.

\bibitem{hardy1917average}
GH~Hardy.
\newblock The average order of the arithmetical functions p (x) and $\delta$
  (x).
\newblock {\em Proceedings of the London Mathematical Society}, 2(1):192--213,
  1917.

\bibitem{heath1992distribution}
DR~Heath-Brown.
\newblock The distribution and moments of the error term in the dirichlet
  divisor problem.
\newblock 1992.

\bibitem{heath2012lattice}
DR~Heath-Brown.
\newblock Lattice points in the sphere.
\newblock In {\em Number theory in progress}, pages 883--892. de Gruyter, 2012.

\bibitem{huxley2003exponential}
Martin~N Huxley.
\newblock Exponential sums and lattice points iii.
\newblock {\em Proceedings of the London Mathematical Society}, 87(3):591--609,
  2003.

\bibitem{iwaniec1988divisor}
Henryk Iwaniec and CJ~Mozzochi.
\newblock On the divisor and circle problems.
\newblock {\em Journal of Number theory}, 29(1):60--93, 1988.

\bibitem{Skriganov}
M.~Skriganov.
\newblock Ergodic theory on sl(n), diophantine approximations and anomalies in
  the lattice point problem.
\newblock {\em Inventiones Mathematicae}, 132:1--72, 04 1998.

\bibitem{skriganov1994constructions}
Maxim~Mikhailovich Skriganov.
\newblock Constructions of uniform distributions in terms of geometry of
  numbers.
\newblock {\em {\selectlanguage{russian} Алгебра и анализ}},
  6(3):200--230, 1994.

\bibitem{sprindzhuk1979metric}
Vladimir~Gennadievich Sprindzhuk.
\newblock {\em Metric theory of Diophantine approximations}.
\newblock VH Winston, 1979.

\bibitem{trevisan2021}
Julien Trevisan.
\newblock Lattice counting problem, 2021.

\bibitem{trevisan2021limit}
Julien Trevisan.
\newblock Limit laws in the lattice problem. ii. the case of ovals, 2021.

\end{thebibliography}
\end{document}